\documentclass{article}
\usepackage{graphicx,amsfonts,amsmath,mathrsfs,amssymb,amsthm,url,color}
\usepackage{fancyhdr,indentfirst,bm,enumerate,natbib}
\usepackage{tikz}
\usepackage[colorlinks=true,citecolor=blue]{hyperref}
\usepackage{dsfont}
\def\Fbar{ {\overline F}}         
          \def\Gbar{ {\overline G}}
\def\oo{\infty}                   \def\d{\,\mathrm{d}}
\def\lm{\lambda}                  \def\th{\theta}
                 
\def\dfrac{\displaystyle \frac}
\newcommand{\id}{\mathds{1}}

\newcommand{\E}{\mathbb{E}}
\newcommand{\R}{\mathbb{R}}
\newcommand{\N}{\mathbb{N}}
\newcommand{\VaR}{\mathrm{VaR}}

\newcommand{\essinf}{\mathrm{ess\mbox{-}inf}}

\renewcommand{\P}{\mathbb{P}}
\renewcommand{\(}{\left (}
\renewcommand{\)}{\right )}
\renewcommand{\[}{\left [}
\renewcommand{\]}{\right ]}

\newtheorem{theorem}{Theorem}[section]

\newtheorem{lemma}[theorem]{Lemma}
\newtheorem{proposition}[theorem]{Proposition}
\newtheorem{definition}{Definition}[section]

\newtheorem{example}[theorem]{Example}

\newtheorem{remark}[theorem]{Remark}

\numberwithin{equation}{section}
\numberwithin{theorem}{section}

\renewcommand{\cite}{\citet}
\allowdisplaybreaks

\def\red{\color{red}}

\begin{document}
	
\baselineskip=18pt
	
\title{Further Developments on Stochastic Dominance for Convex\\ Combinations of Infinite-Mean Random Variables}
	
\author{ Keyi Zeng\thanks{School of Management, University of Science and Technology of China, China. 
Email: \url{kyzeng@mail.ustc.edu.cn}} 
\and Zhenfeng Zou\thanks{School of Public Affairs, University of Science and Technology of China,  China. Email: \url{zfzou@ustc.edu.cn}}
\and Yuting Su\thanks{School of Management, University of Science and Technology of China,  China.
Email: \url{syt20020224@mail.ustc.edu.cn}}
\and Taizhong Hu\thanks{School of Management, University of Science and Technology of China,  China. 
Email: \url{thu@ustc.edu.cn}}
  	}

\date{November, 2025 \\ Revised April, 2026}
	
\maketitle

\begin{abstract}

In recent years, stochastic dominance for independent and identically distributed (iid) infinite-mean random variables has received considerable attention. The literature has identified several classes of distributions of nonnegative random variables that encompass many common heavy-tailed distributions. A key result demonstrates that the weighted sum of iid random variables from these classes is stochastically larger than any individual random variable in the sense of the first-order stochastic dominance. This paper systematically investigates the properties and inclusion relationships among these distribution classes, and extends some existing results to more practical scenarios. Furthermore, we analyze the case where each random variable follows a compound binomial distribution, establishing necessary and sufficient conditions for the preservation of the aforementioned stochastic dominance relation.

\medskip
		
 		 \medskip
		
		
\end{abstract}

\section{Introduction}

The study of stochastic comparisons between linear combinations of random variables is a cornerstone of probability theory with profound implications in fields such as finance, insurance, and economics. A vast literature has explored this topic, typically under the assumption that the random variables involved possess finite expectations \citep[see, e.g.][and the references therein]{Pro65, Ma00, XH11}.

However, a growing body of research has revealed intriguing and counter-intuitive phenomena when the finite expectation is relaxed. In a seminal contribution, \cite{Ibr05} showed that for independent and identically distributed (iid) random variables $X_1, \ldots, X_n$ with infinite-mean stable distributions on the positive side, a more diversified portfolio can actually be stochastically larger than a less diversified one. Specifically, he demonstrated that
$$
   \(\sum_{i=1}^n\theta_i\) X_1  \le_{\rm st} \sum_{i=1}^n\theta_i X_i   \eqno {\rm (SD^*)}
$$
for any nonnegative real vector $\bm\theta=(\th_1, \ldots, \th_n)$, where $\le_{\rm st}$ is the \emph{usual stochastic order} or the \emph{first-order stochastic dominance}. For iid random variables $X_1, \ldots, X_n$ having s symmetric stable distribution with infinite mean, \cite{Ibr09} established 
$$
     \(\sum_{i=1}^n\theta_i\) |X_1|  \le_{\rm st} \left |\sum_{i=1}^n\theta_i X_i\right |
$$
for any nonnegative vector $\bm\theta$. This surprising result ${\rm ({\red SD^*})}$, which runs counter to the risk-diversification paradigm of traditional finance, was later extended to other important heavy-tailed models, including the Pareto distribution \citep{EMS02, CEW25}. In view of Value-at-Risk, ${\rm ({\red SD^*})}$ implies that independence is worse than perfect dependence no matter how large we choose the confidence level. The definition of the Pareto distribution is given in Section \ref{sect-2}. In portfolio diversification, property ${\rm ({\red SD^*})}$ has an intuitive implication: A more diversification portfolio is stochastically larger. 

This intriguing finding spurred a line of research aimed at identifying the broadest possible distributional families for which such stochastic dominance relations hold. A key contribution was made by \cite{CHSZ25}, who introduced the class $\mathcal{H}$, characterized by the concavity of a transformed survival function, and the broader class $\mathcal{H}^\ast$ (also studied by \cite{ALO24}), defined via subadditivity. Subsequently, \cite{CS25} proposed the class $\mathcal{G}$, based on the subadditivity of $-\log F(1/x)$, offering a different perspective on the phenomenon. In a unifying effort, \cite{Mul25} demonstrated that many of these results, including those in \cite{CEW25} and \cite{CS25}, can be seen as special cases of a broader ``super-Cauchy" framework, thereby clarifying the relationships between super-Pareto, super-Fr{\'e}chet, and super-Cauchy families. Independently, \cite{Vin25} developed the class $\mathcal{V}$, defined by a ``completely subscalable" property of the survival function $\overline{F}(x)$, to address a distinct yet related dominance result for a ``concentrated portfolio" model. 
(The formal definitions of the classes $\mathcal{V}$, $\mathcal{H}$, $\mathcal{H}^\ast$, and $\mathcal{G}$ are provided in Section \ref{sect-2}, while the super-Pareto, super-Fr{\'e}chet, and super-Cauchy families are defined in Remark \ref{re-20251031}.)

While these classes are known to encompass many common heavy-tailed distributions, the literature has so far lacked a systematic and unified treatment of the relationships between them. Scattered results suggest inclusions---for instance, it is known that $\mathcal{G} \subset \mathcal{H}^\ast$ \citep{ALO24} and that $\mathcal{H} \subset \mathcal{V}$---but a clear, comprehensive picture of their hierarchy and relative sizes is missing. 
For example, it remains to be seen whether $\mathcal{V}$ is contained in $\mathcal{H}^\ast$, or how $\mathcal{G}$ intersects with $\mathcal{H}$ and $\mathcal{V}$. Establishing this hierarchy is not merely a taxonomic exercise; it is essential for both theoretical clarity and practical model selection. Understanding, for instance, that a result proven for the broader class $\mathcal{V}$ automatically applies to the narrower class $\mathcal{H}$ allows researchers to build upon the most general findings.

This paper aims to fill this gap by providing a comprehensive analysis of these four classes of infinite-mean distributions. Our main contributions are fourfold.
\begin{itemize}
 \item \textbf{First, we systematically map the inclusion relationships among the classes $\mathcal{H}$, $\mathcal{V}$, $\mathcal{H}^*$, and $\mathcal{G}$} (Propositions \ref{pr-250925} and \ref{pr-250928}).
We provide rigorous proofs of previously implicit relationships, such as $\mathcal{V} \subsetneq \mathcal{H}^\ast$, and construct concrete counterexamples to delineate their boundaries. This culminates in a detailed Venn diagram that precisely illustrates the hierarchy, clarifying, for instance, the proper inclusions $\mathcal{H} \subsetneq \mathcal{V} \subsetneq \mathcal{H}^\ast$ and the non-trivial intersection of $\mathcal{G}$ with the others. 
We also investigate the closure under power transformations of distribution and survival functions, maximum transformations of random variables,  convex transformations of random variables, and others. Some of these properties are already known, while others are newly established.

 \item \textbf{Second, we extend the classical stochastic dominance result Theorem \ref{th-H} to more realistic and practical scenarios}. 
Building on the framework of \cite{CHWZ25}, we consider settings where losses are triggered by rare events (Theorems \ref{th-251002} and \ref{th-251003}), where variables have heavy tails only above a threshold (Proposition \ref{pr-251002}), and where risks are subject to an upper bound (Proposition \ref{pr-250930}). These extensions significantly enhance the applicability of the theory.

 \item \textbf{Third, we provide two complementary results that sharpen the understanding of the existing theorems}. 
First, we present a counterexample demonstrating that the class $\mathcal{H}$ in Theorem \ref{th-H} cannot be replaced by the larger class $\mathcal{V}$, thereby establishing the optimality of the conditions in the original result. Second, we offer a simple and concise proof of Theorem \ref{th-V} using an induction method, which not only clarifies the underlying mechanism but also serves as a useful template for establishing similar dominance relations in other settings.

 \item \textbf{Finally, we investigate the preservation of stochastic dominance relations in the context of compound distributions}. 
While \cite{CHSZ25} established results for compound Poisson models, we extend this line of inquiry to compound binomial distributions, providing necessary and sufficient conditions (Theorem \ref{th-251001}) that mirror and generalize those for the Poisson case.
\end{itemize}

The rest of the paper is organized as follows. 
Section \ref{sect-2} provides the necessary preliminaries, including the formal definitions of the four distribution classes $\mathcal{H}$, $\mathcal{V}$, $\mathcal{H}^\ast$, and $\mathcal{G}$, the concept of majorization order, and a review of the known stochastic dominance results (Theorems \ref{th-H}–\ref{th-V}) that serve as the foundation for our work. In Section \ref{sect-3}, we systematically investigate the properties of these distribution classes and establish the precise inclusion relationships among them. Section \ref{sect-4} presents our main extensions of the stochastic dominance results to more practical settings, including losses triggered by rare events, losses with heavy tails only above a threshold, and truncated random variables. Section \ref{sect-5} examines the case of compound binomial distributions, providing necessary and sufficient conditions for the preservation of the ({\red SD}) and ${\rm ({\red SD^*})}$ relations. Finally, Section \ref{sect-6} concludes the paper with a discussion of open problems and directions for future research. Some technical proofs of the results in Section \ref{sect-3} are relegated to the Appendix.

Throughout, random variables are defined on an atomless probability space $(\Omega, \mathscr{F}, \P)$. We write $\bm X\stackrel {d}{=} \bm Y$ if $\bm X$ and $\bm Y$ have the same distribution, and write $f(x)\stackrel {\rm sgn}{=} g(x)$ if two functions $f(x)$ and $g(x)$ have the same sign. For a distribution function $F$, its left-continuous inverse is defined by
$$
   F^{-1}(\alpha) =\inf \{x\in \R: F_X(x)\ge \alpha,\quad \alpha\in (0,1],
$$
with $F^{-1}(0)=\inf\{x\in\R: F(x)>0\}$.  Denote by $\N$ the set of all positive integers, $\R_+$ be the set of all nonnegative real number, and $\R_{++}$ be the set of all positive real numbers. For $n\in\N$, let $[n]=\{1, \ldots, n\}$. Denote $\Delta_n=\left\{\bm\th\in (0,1)^n: \sum^n_{i=1} \th_i=1\right\}$. Also, ``increasing'' and ``decreasing'' mean ``nondecreasing'' and ``nonincreasing'', respectively. The ratio $a/0$ is understood to be $+\oo$ whenever $a>0$, and the ratio $0/0$ is not well-defined.

\section{Preliminaries and known results}
\label{sect-2}

\subsection{Definitions}
First, we introduce some concepts and terminology to be used in the sequel. A function $\varphi$ is said to be subadditive if $\varphi(x+y)\le\varphi(x) + \varphi(y)$ for all $x, y$ in the domain of $\varphi$. The function $\varphi$ is said to be superadditive if the inequality is reversed. A function $\varphi: \R_+\to\R$ is said to be star-shaped if $\varphi(0) \le 0$ and $\varphi(x)/x$ is increasing in $x\in\R_{++}$. If $\varphi(0) \ge 0$ and $\varphi(x)/x$ is decreasing in $x\in\R_{++}$, then $\varphi$ is said to be anti-star-shaped.

The notion of majorization defines a partial ordering of the diversity of the components of vectors. To recall the definition of majorization order \citep{MOA11}, let $a_{(1)}\le a_{(2)}\le\cdots \le a_{(n)}$ be the
increasing arrangement of components of the vector $\bm a=(a_1, a_2,\ldots, a_n)$. For vectors $\bm a, \bm b \in \R^n$,  $\bm a$ is said to be majorized by $\bm b$, denoted by $\bm a\preceq_{\rm m} \bm b$, if $\sum_{i=1}^n a_{(i)}=\sum_{i=1}^n b_{(i)}$ and
\begin{equation}
\label{eq-250912}
     \sum_{i=1}^j a_{(i)}\ge \sum_{i=1}^j b_{(i)}\ {\rm for}\ j\in [n-1].
\end{equation}
If the strict inequality \eqref{eq-250912} holds for at least one $j\in [n-1]$,  $\bm a\preceq_{\rm m} \bm b$ is denoted by $\bm a\prec_{\rm m} \bm b$. A real-valued function $\phi$ defined on a set $A\subseteq \R^n$ is said to be Schur-concave [Schur-convex] on $A$ if $\phi(\bm a) \ge [\le]\, \phi(\bm b)$ whenever $\bm a\preceq_{\rm m} \bm b$ and $\bm a, \bm b\in A$.

Throughout this paper, we always assume random variables are nonnegative unless stated otherwise.

\begin{definition}
{\rm \citep{Vin25}}\ \ Let $F$ be a non-degenerate distribution function with $F(0-)=0$. $F$ is said to be completely subscalable if the inequality
\begin{equation}
 \label{eq-250903}
      \th\, \Fbar(x) \le \Fbar\(\frac {x}{\th}\)
\end{equation}
holds for all $x\in\R_+$ and all $\th\in (0,1)$. Denote by ${\mathcal V}$ the class of all completely subscalable  distribution functions. 
\end{definition}

The property \eqref{eq-250903} is equivalent to the \emph{anti-starshapedness} of $\Fbar(1/x)$ on $\R_{++}$, i.e.,
\begin{equation*}
   \label{quasi-homo}
   \Fbar\(\frac{x}{t}\) \leq t \Fbar\(x\),\quad x \in \R_+,\ t > 1.
\end{equation*}
For more details on anti-starshapedness, we refer the reader to \cite{Lynch87}. 

\begin{definition}\label{def-H}
{\rm \citep{CHSZ25}}\ \ Let $F$ be a non-degenerate distribution function with $F(0-)=0$. We say $F\in \mathcal{H}$ if the function $\Fbar (1/x)$ is concave in $x \in \R_{++}$.
\end{definition}

\begin{definition}
{\rm \citep{ALO24, CHSZ25}}\ \ Let $F$ be a non-degenerate distribution function with $F(0-)=0$. We say $F\in \mathcal{H}^\ast$ if $\Fbar (1/x)$ is subadditive in $x\in \R_{++}$. 
\end{definition}

\begin{definition}
{\rm \citep{CS25}}\ \ Let $F$ be a non-degenerate distribution function with $F(0-) = 0$. We say $F \in \mathcal{G}$ if the function
\begin{equation}
 \label{eq-250927}
   \Lambda_F(x) = -\log F\(\frac{1}{x}\)
\end{equation}
is subadditive in $x \in \R_{++}$ with the convention $\log 0=-\oo$. 
\end{definition}

The distributions in $\mathcal{H}^\ast$ are called InvSub (inverted subadditivity) by \cite{ALO24} who also showed that $\mathcal{H}^\ast$ is more general than the class of super-Pareto distributions. Analogously, \cite{CHSZ25} termed the distributions in $\mathcal H$ as InvCv (inverted concavity). Clearly, $\mathcal{H} \subset \mathcal{H}^\ast$.  In fact, $\mathcal{H}$ is a proper subset of $\mathcal{H}^\ast$ (see Example \ref{ex-250920}). Equivalent characterizations of distributions in $\mathcal{V}$ and $\mathcal{H}$ are as follows:
\begin{itemize}
  \item \citep{Vin25}\ \ $F\in \mathcal{V}$ if and only if $x \Fbar(x)$ is increasing in $x\in\R_+$.

  \item \citep[][Proposition 2]{CHSZ25}\ \ Let $F$ have a density function $f$, then $F\in \mathcal{H}$ if and only if $x^2f(x)$ is increasing in $x\in\R_+$.

  \item \citep[][Proposition 2.5]{ALO24}\ \ Let $F$ have a density function $f$, and $\lm(t)=f(t)/\Fbar(t)$ denote the failure rate of $F$. If $ x \lm(x) \le 1$ for all $x\in \R_+$, then $F\in \mathcal{H}^\ast$.
\end{itemize}

For nonnegative functions on $\mathbb{R}_{+}$, it is well known that increasing concave functions are anti-starshaped, and increasing anti-starshaped functions are subadditive \citep[see, e.g.][Chapter 4.B]{SS07}. 
Consequently, we obtain the hierarchical inclusion relationships 
$$
    \mathcal{H} \subset \mathcal{V} \subset \mathcal{H}^\ast.
$$
While these inclusions follow from the general functional hierarchy, we also provide direct proofs based on the definitions of the respective classes in Section \ref{sect-3}, where the strictness of these inclusions is further demonstrated through concrete counterexamples. This functional perspective not only clarifies the structure of these classes but also provides a unifying framework for understanding many of their properties.

Regarding the class $\mathcal{G}$, its definition via the subadditivity of $-\log F(1/x)$ places it in a different but related hierarchy. \cite[Section 4]{ALO24} discuss the relationship between $\mathcal{G}$ and $\mathcal{H}^\ast$, showing that $\mathcal{G} \subset \mathcal{H}^\ast$ and that the inclusion is strict.

For $\alpha>0$, the Pareto distribution, denoted by Pareto$(\alpha)$, is given by
$$
    F_\alpha(x) =1- \frac {1}{x^\alpha},\quad x\ge 1,
$$
and  the Fr\'{e}chet distribution, denoted by Fr\'{e}chet$(\alpha)$, is given by
$$
   F_\alpha(x) =\exp\left\{- x^{-\alpha}\right\},\quad x>0.
$$
For $\alpha\in (0, 1]$, both distributions have infinite means, and belong to any one of $\mathcal{H}$ and $\mathcal{G}$. Many other examples of distributions in $\mathcal{H}$ and $\mathcal{G}$ are listed in \cite{CHSZ25} and \cite{CS25}, respectively.

\begin{remark}\label{continue}
{\rm (\emph{Continuity of $F$ on $\R_{++}$})\ \ 
\cite[Lemma 5.4]{Vin25} proved that $F\in \mathcal{V}$ implies $F(x)$ is continuous on $\R_{++}$. 
A similar continuity property holds for the class $\mathcal H$. Indeed, if $F\in\mathcal{H}$, then by Definition \ref{def-H}, the function $\eta(x) := \Fbar(1/x)$ is concave on $\R_{++}$.  A concave function on an open interval is necessarily continuous \citep[see, e.g.,][Theorem 10.1]{Rock70}. Consequently, $\eta(x)$ is continuous on $\R_{++}$, which immediately implies that $F(x)$ is continuous on $\R_{++}$ as well. Thus, all distributions in $\mathcal{H}$ are continuous on $\R_{++}$.  }
\end{remark}

\begin{remark}[{\it Essential infimum}]
{\rm Note that $F\in\mathcal{G}$ is equivalent to
\begin{equation}
  \label{eq-250906}
    F\(\frac {xy}{x+y}\) \ge F(x)  F(y),\quad (x,y)\in\R_{++}^2.
\end{equation}
This implies $\essinf (F)=0$, that is, $F(x)>0$ for any $x\in\R_{++}$. Now let $X$ be a truncated Fr\'{e}chet random variable with density function given by
$$
   f(x) = \left\{\begin{array}{ll} 0, & x\in [0,1], \\
      c x^{-\alpha-1} \exp\{-x^{-\alpha}\}, & x> 1,  \end{array}\right.
$$
where $c>0$ is a normalized constant. Then $F\in \mathcal{H}$ since $x^2 f(x)$ is increasing in $x\in\R_+$. Thus, in view of Proposition \ref{pr-250925} (i), $F\in \mathcal{H}$, $\mathcal{V}$ or $\mathcal{H}^\ast$ does not necessarily imply $\essinf(F)=0$. Example \ref{ex-250922} also shows that $F\in\mathcal{V}$ or $F\in\mathcal{H}^\ast$ does not necessarily imply $\essinf(F)=0$. In view of these observations, we have $\mathcal{H} \not \subset \mathcal{G}$.   }
\end{remark}

\begin{remark}[{\it Transform order characterizations}]
\label{convex_trans}
{\rm Let $X$ and $Y$ be two nonnegative random variables with distribution functions $F$ and $G$, respectively.
We say that $X$ is smaller than $Y$ 
\begin{itemize}
\item[(1)] in the convex transform order, denoted by $X \le_{\rm c} Y$ or $F \le_{\rm c} G$, if $G^{-1}(F(x))$ is convex in $x \in \R_+$;
    
\item[(2)] in the star order, denoted by $X \le_* Y$ or $F \le_* G$, if $G^{-1}(F(x))/x$ is increasing in $x \in \mathbb{R}_{++}$;

\item[(3)] in the superadditive order, denoted by $X \le_{\rm su} Y$ or $F \le_{\rm su} G$, if $G^{-1}(F(x))$ is superadditive on $\R_+$.
\end{itemize}
For a comprehensive treatment of these stochastic orders, we refer the reader to \cite{SS07}. The order $F\le_{\rm c} G$ gives us an intuition that $F$ is less skewed to the right than $G$. This concept is discussed in detail in \cite{Zwe64} and \cite{BP81}.

Now, the classes $\mathcal{V}$, $\mathcal{H}$, and $\mathcal{H}^\ast$ admit alternative characterizations in terms of stochastic transform orders. Recall that for a distribution function $F$ with $F(0-)=0$, the random variable $1/X$ (where $X \sim F$) plays a key role. Specifically, we have the following equivalences: 
\begin{align*}
   F \in \mathcal{H} & \iff U \leq_{\rm c} \frac {1}{X}, \\[3pt]
   F \in \mathcal{V} & \iff U \leq_* \frac {1}{X}, \\[3pt]
   F \in \mathcal{H}^\ast & \iff U \leq_{\rm su} \frac {1}{X},
\end{align*}
where $U$ is the uniform distribution on $(0,1)$.  }
\end{remark}

In the sequel, a random variable $X$ is said to be $\mathcal{T}$-distributed if its distribution function belongs to the class $\mathcal{T}$, where $\mathcal{T}$ can be any one of $\mathcal{H}$, $\mathcal{V}$, $\mathcal{H}^\ast$ or $\mathcal{G}$. For $X \sim F$ where $F \in\mathcal{T}$, we also write $X \in \mathcal{T}$.

\subsection{Known results on stochastic dominance}

Several fundamental dominance results have been established for the distribution classes introduced above. 
These theorems, which serve as the primary motivation and foundation for our work, are summarized below.

The first result, due to \cite{CHSZ25}, establishes stochastic dominance for linear combinations of iid random variables from the classes $\mathcal{H}$.
\begin{theorem}
  \label{th-H}
{\rm \citep{CHSZ25}}\ \ Let $\bm X=(X_1, \ldots, X_n)$ be a vector of iid random variables with a common distribution function $F$. If $F\in\mathcal{H}$, then ${\rm ({\red SD})}$ holds, i.e.,
$$
   \sum_{i=1}^n\eta_i X_i  \le_{\rm st} \sum_{i=1}^n\theta_i X_i  \eqno {\rm (SD)}
$$
for $\bm \theta, \bm \eta\in \R_+^n$ such that $\bm\theta \preceq_{\rm m} \bm \eta$. 
\end{theorem}

A stronger dominance relation, known as ${\rm ({\red SD^*})}$, holds for the larger class $\mathcal{H}^\ast$. 
This result was established independently by \cite{ALO24} and \cite{CHSZ25}.

\begin{theorem}
  \label{th-Hast}
{\rm \citep{ALO24, CHSZ25}}\ \ Let $\bm X=(X_1, \ldots, X_n)$ be a vector of iid random variables with a common distribution function $F$.
If $F\in\mathcal{H}^\ast$, then ${\rm ({\red SD^*})}$ holds for all $\bm\theta\in \R_+^n$, i.e.,
$$
   \(\sum_{i=1}^n\th_i\) X_1 \le_{\rm st} \sum_{i=1}^n\theta_i X_i  \eqno {\rm (SD^*)}
$$
\end{theorem}

The same dominance relation ${\rm ({\red SD^*})}$ also holds for the class $\mathcal{G}$, as shown by \cite{CS25}.

\begin{theorem}
  \label{th-G}
{\rm \citep{CS25}}\ \ Let $\bm X=(X_1, \ldots, X_n)$ be a vector of iid random variables with a common distribution function $F\in\mathcal{G}$. Then ${\rm ({\red SD^*})}$ holds for all $\bm\theta\in \R_+^n$.
\end{theorem}

A different type of dominance result, concerning a ``concentrated portfolio" model, was established by \cite{Vin25} for the class $\mathcal{V}$.

\begin{theorem}
 \label{th-V}
{\rm \citep{Vin25}}\ \ Let $\bm X=(X_1, \ldots, X_n)$ be a vector of independent random variables with $X_i\sim F_i\in\mathcal{V}$ for each $i$, and let $(I_1, \ldots, I_n)$ be a multivariate Bernoulli random vector, independent of $\bm X$, satisfying $\sum^n_{i=1} I_i=1$ and $\P(I_i=1)=\theta_i$ for each $i$, where $\sum^n_{i=1} \th_i=1$. Then
$$
   \sum^n_{i=1} I_i X_i \le_{\rm st} \sum^n_{i=1} \th_i X_i.    \eqno {\rm (SD}_{\rm cp})
$$
\end{theorem} 

\begin{remark}
{\rm 
In Theorem \ref{th-V}, the vector $\bm I$ has exactly one component equal to $1$ and the others equal to $0$. The sum $\sum^n_{i=1} I_i X_i$ is termed as \emph{concentrated portfolio} by \cite{Vin25}, which concentrates all exposure on a single risk (i.e., selects exactly one of the $X_i$ at random according to the random weights). Thus, the stochastic dominance between the diversified portfolio $\sum^n_{i=1} \th_i X_i$ and the concentrated portfolio $\sum^n_{i=1} I_i X_i$ is referred to property ${\rm ({\red SD}}_{\rm \red cp})$.  If $X_1,\ldots, X_n$ are iid, then $ {\rm ({\red SD}_{\rm \red cp}})$ reduces to $({\red {\rm SD}^\ast})$.

Actually, \cite{Vin25} established a more general result than Theorem \ref{th-V}, which is called the one-basket-theorem. This theorem states as follows:
Let $\bm X=(X_1, \ldots, X_n)$ be a vector of independent random variables with $X_i\sim F_i$ for each $i$.
Given a weight vector $\bm \th \in \Delta_n$, let $(I_1, \ldots, I_n)$ be defined as in Theorem \ref{th-V}. 
Suppose that for each $i \in [n]$ and every $B\subset [n]$ with $B\supseteq \{i\}$, $F_i$ satisfies  
\begin{equation*}
\theta_B \overline{F}_i(x) \leq \overline{F}_i(x / \theta_B),\quad x\in\R_+, 
\end{equation*}
where $\th_B = \sum_{j \in B} \th_j$ is the subset weight of $\mu$. Then (${\rm SD}_{\rm cp}$) holds.    
}
\end{remark}

\section{Properties of distribution classes}
\label{sect-3}

If $X$ belongs to any of the classes $\mathcal{H}$, $\mathcal{V}$, $\mathcal{H}^\ast$ and $\mathcal{G}$, then $cX$ also belongs to the same class for $c\in\R_{++}$. Further properties of these four classes are listed in the following four propositions (Propositions \ref{pr-251001},  \ref{pr-251001x}, \ref{pr-250925}, and \ref{pr-250928}).

Many commonly encountered examples, such as those listed in Table 1 of \cite{CS25}, satisfy the condition that $\Lambda_F(x)$ is a concave function on $\R_{++}$.

\begin{proposition}
 \label{pr-251001}
If $\Lambda_F(x)$ is a concave function on $\R_{++}$, then $F^\beta \in \mathcal{H}$ [resp. $\mathcal{V}$ and $\mathcal{H}^\ast]$ for all $\beta \in (0,1)$.
\end{proposition}

\begin{proof}
Define $\eta(x)= \Fbar(1/x)$. It suffices to show that $\eta_\beta(x)= 1-\[1 - \eta(x)\]^\beta$ is concave. Observing that $\Lambda_F(x) = - \log F(1/x) = -\log\(1 - \eta(x)\)$, we have $\eta_\beta(x) = 1 - \exp\{-\beta \Lambda_F(x)\}$. Since the function $t \mapsto 1 - \exp(-\beta t)$ is increasing and concave, and $\Lambda_F$ is concave by assumption, it follows that $\eta_\beta(x)$ is concave as a composition of a concave and  increasing function with a concave function.
\end{proof}

\begin{proposition}
 \label{pr-251001x}
Let $X$ be an absolutely continuous random variable. If $X$ has a decreasing failure rate $\lm(t)$ on $\R_+$, then $1/X\in \mathcal{G}$.
\end{proposition}

\begin{proof}
Denote by $F$ and $G$ [resp. $f$ and $g$] the distribution [resp. density] functions of $X$ and $1/X$, respectively. Then $G(x)=\Fbar (1/x)$, and $g(x)=f(1/x)/x^2$ for $x>0$. To prove $1/X\in\mathcal{G}$, it suffices to verify that $-\log G(1/x)$ is anti-star-shaped, that is, the function $\phi(x)=x\log G(x)$ is decreasing on $\R_+$. Note that, for any $x>0$,
\begin{align*}
   \phi'(x) & =\log G(x) + \frac {xg(x)}{G(x)} =\log\Fbar\(\frac {1}{x}\) +\frac {1}{x} \frac {f(1/x)}{\Fbar(1/x)} \\
   & = - \int^{1/x}_0 \lm(u)\d u + \frac {1}{x} \lambda \(\frac {1}{x}\)\\
    & =  \int^{1/x}_0 \[ \lambda \(\frac {1}{x}\) -\lm(u)\] \d u\le 0,
\end{align*}
where the inequality follows since $\lm(t)$ is decreasing. Therefore, the desired result follows.  
\end{proof}

It is well-known \citep[see, e.g.,][]{BP81} that if $X$ has a log-convex density function on $\R_+$, then $X$ has a decreasing failure rate. Proposition \ref{pr-251001x} provides a sufficient condition for verifying whether a distribution belongs to $\mathcal{G}$. For example, let $X$ follow $\Gamma (\alpha, \beta)$ distribution with shape parameter $\alpha\in \R_{++}$ and scale parameter $\beta\in\R_{++}$. Then $1/X$ follows the inverse-$\Gamma(\alpha, \beta)$ distribution with density function
$$
   g(x)= \frac {\beta^\alpha x^{-\alpha-1}}{\Gamma(\alpha)} \exp\left\{-\frac {\beta}{x}\right \},\quad x\in\R_{++}.
$$
Since $X$ has a log-convex density function when $\alpha\in (0,1]$, by Proposition \ref{pr-251001x}, $1/X\in \mathcal{G}$. 

To state the next proposition, we recall from \cite{SS07} the definitions of some common used stochastic orders.
For two random variables $X$ and $Y$ with respective distribution functions $F_X$ and $F_Y$, $X$ is said to be smaller than $Y$ in the \emph{hazard rate order}, denoted by $X\le_{\rm hr} Y$ or $F_X\le_{\rm hr} F_Y$, if $\Fbar_Y(x)/\Fbar_X(x)$ is increasing in $x$ for which the ratio is well-defined. $X$ is said to be smaller than $Y$ in the \emph{likelihood ratio order}, denoted by $X\le_{\rm lr} Y$ or $F_X\le_{\rm lr} F_Y$, if $F_X$ and $F_Y$ have the density functions $f_X$ and $f_Y$, respectively, satisfying that $f_Y(x)/f_X(x)$ is increasing in $x$ for which the ratio is well-defined. 

\begin{proposition}
\label{pr-250925} \
\begin{itemize}
  \item[{\rm (i)}] $\mathcal{G} \varsubsetneq \mathcal{H}^\ast$ \citep[][Theorem 4.13]{ALO24}.
      $\mathcal{H}\varsubsetneq \mathcal{V} \varsubsetneq  \mathcal{H}^\ast$.

  \item[{\rm (ii)}]  If $F\in\mathcal{H}$, then $F^\beta\in\mathcal{H}$ for $\beta\ge 1$
     \citep[][Proposition 3 (i)]{CHSZ25}. \\
     If $F\in \mathcal{G}$, then $F^\beta \in \mathcal{G}$ for all $\beta>0$ \citep[][Proposition 2(ii)]{CS25}. \\
      If $F\in\mathcal{V}$  [resp. $\mathcal{H}^\ast$], then $F^\beta\in \mathcal{V}$  [resp. $\mathcal{H}^\ast$] for all $\beta\ge 1$.

  \item[{\rm (iii)}] If $F\in\mathcal{V}$ [resp. $\mathcal{H}$, $\mathcal{H}^\ast$, $\mathcal{G}$], then $1-\Fbar^\beta\in \mathcal{V}$ [resp. $\mathcal{H}$, $\mathcal{H}^\ast$, $\mathcal{G}$] for all $\beta\in (0,1)$.

  \item[{\rm (iv)}]  If $F\in\mathcal{H}$ and $F\le_{\rm lr} G$, then $G \in \mathcal{H}$
       \citep[][Proposition 3(iv)]{CHSZ25}.\\
      If $F\in\mathcal{V}$ [resp. $\mathcal{H}^\ast$] and $F\le_{\rm hr} G$, then $G \in \mathcal{V}$ [resp. $\mathcal{H}^\ast$].

  \item[{\rm (v)}]  For $w_1, \ldots, w_n\in\R_+$ such that $\sum^n_{i=1} w_i=1$,
    \begin{itemize}
      \item if $F_1, \ldots, F_n\in \mathcal{V}$, then $\sum^n_{i=1} w_i F_i \in \mathcal{V}$.
      \item if $F_1, \ldots, F_n\in \mathcal{H}$ [resp. $\mathcal{H}^\ast$], then $\sum^n_{i=1} w_i F_i \in \mathcal{H}$ [resp. $\mathcal{H}^\ast$] \citep[][Proposition 4]{CHSZ25}.
      \item If $F_1, \ldots, F_n \in \mathcal{G}$ and $F_1\le_{\rm st} \cdots \le_{\rm st} F_n$, then $\sum_{i=1}^n w_i F_i \in \mathcal{G}$ \citep[][Proposition 3]{CS25}.
    \end{itemize}
\end{itemize}
\end{proposition}

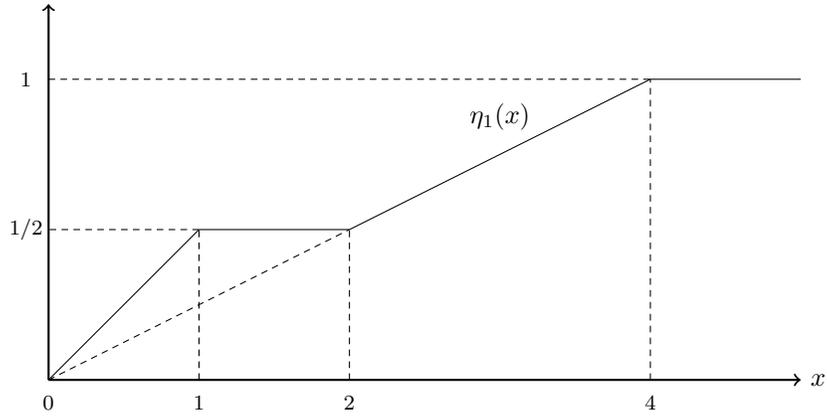
\begin{figure}[htbp]
	\centering
	\begin{minipage}[b]{0.5\textwidth}
		\begin{tikzpicture}
			\tikzset{ box/.style ={ circle,
					minimum width =1.5pt,
					minimum height =1.5pt,
					inner sep=1.5pt,
					draw=black,
					fill=white   	}    }
			\draw[->,thick] (0,0) -- (10,0) node[right] {$x$};
			\draw[->,thick] (0,0) -- (0,5) node[above] { };	
			\draw[domain=0:2,color=black] plot (\x,\x);
			\draw[domain=2:4,color=black] plot (\x,2);
			\draw[domain=4:8,color=black] plot (\x, {2+0.5*(\x-4)});
            \draw[domain=8:10,color=black] plot (\x,4);
			\node at (6,3.5){$\eta_1(x)$};
			\draw  [dash pattern={on 2.5pt off 2pt}] (0,0) -- (4,2);
     		\draw  [dash pattern={on 2.5pt off 2pt}] (0,4) -- (8,4);
			\draw  [dash pattern={on 2.5pt off 2pt}] (2,2) -- (0,2);
			\draw  [dash pattern={on 2.5pt off 2pt}] (2,0) -- (2,2);
			\draw  [dash pattern={on 2.5pt off 2pt}] (4,0) -- (4,2);
            \draw  [dash pattern={on 2.5pt off 2pt}] (8,0) -- (8,4);
			\node at (2,-0.3) {\footnotesize $1$};
			\node at (4,-0.3) {\footnotesize $2$};
			\node at (8,-0.3) {\footnotesize $4$};
			\node at (0,-0.3) {\footnotesize $0$};
			\node at (-0.3,2) {\footnotesize $1/2$};
            \node at (-0.3,4) {\footnotesize $1$};
		\end{tikzpicture}
	\caption{The function $\eta_1(x)$}
			\label{Fig-250901}
	\end{minipage}
\end{figure}

\begin{example}
 \label{ex-250920} {\rm
($\mathcal{H} \varsubsetneq \mathcal{V}$ and $\mathcal{V}\not\subset\mathcal{G}$).\ \ Let $F_1$ be a distribution function with $F_1(0-)=0$, and $\eta_1(x)=\Fbar_1(1/x)$ be defined as follows (see Figure \ref{Fig-250901})
$$
   \eta_1(x) =\left\{ \begin{array}{ll} x/2, & x\in [0, 1],\\
       1/2, & x\in (1,2],\\
       x/4, & x\in (2, 4), \\
       1, & x\in [4, \oo).  \end{array}
   \right.
$$
It is easy to see that $\eta_1(x)$ is not concave, and $\eta_1(x)/x = (1/x) \Fbar_1(1/x)$ is decreasing in $x\in\R_+$. Thus, $F_1\not\in\mathcal{H}$, but $F_1\in\mathcal{V}$, implying $\mathcal{H} \varsubsetneq \mathcal{V}$. On the other hand, $F_1\not\in\mathcal{G}$ since $\essinf (F_1)=1/4$.     }
\end{example}

\begin{example}
 \label{ex-250918}
{\rm $(\mathcal{H} \varsubsetneq \mathcal{V})$.\ \ Consider a distribution function $F$ such that $F(x)=0$ for $x<1$, and
\begin{equation*}
   F(x)= 1 - \frac{3(1/x - 1)^2 + 1}{x}, \quad x \ge 1.
\end{equation*}
Denote $g(x)= x \Fbar(x)$. Then $g(x)=x$ for $x\in [0,1]$, and $g(x) = 3(1/x-1)^2 + 1$ for $x>1$. It is easy to see that
\begin{equation*}
    g'(x) = \frac{6}{x^2} \(1 - \frac{1}{x}\)\ge 0,\quad  x \ge 1,
\end{equation*}
implying $g(x)$ is increasing in $x\in (1,\oo)$. Thus, $F \in \cal V$. Denote $\eta(x)=\Fbar(1/x)$. Then $\eta(x) = 3x^3-6x^2+4x$ for $x\in (0,1]$, and $\eta(x)=1$ for $x>1$. Since $\eta''(x) = 6(3x - 2) > 0$ for $x\in (2/3, 1)$, $\eta(x)$ is not concave on $\R_+$, implying $F\not\in\mathcal{H}$. Therefore, $\mathcal{H} \varsubsetneq \mathcal{V}$.  }
\end{example}

\begin{figure}[htbp]
	\centering
	\begin{minipage}[b]{0.5\textwidth}
		\begin{tikzpicture}
			\tikzset{ 	box/.style ={ circle,
					minimum width =1.5pt,
					minimum height =1.5pt,
					inner sep=1.5pt,
					draw=black,
					fill=white   	}    }
			\draw[->,thick] (0,0) -- (10,0) node[right] {$x$};
			\draw[->,thick] (0,0) -- (0,5) node[above] { };	
			\draw[domain=0:2,color=black] plot (\x,\x);
			\draw[domain=2:6,color=black] plot (\x,2);
			\draw[domain=6:8,color=black] plot (\x, {\x-4});
            \draw[domain=8:10,color=black] plot (\x,4);
			\node at (6,3.5){$\eta_2(x)$};
			\draw  [dash pattern={on 2.5pt off 2pt}] (0,0) -- (8,4);
     		\draw  [dash pattern={on 2.5pt off 2pt}] (0,4) -- (8,4);
			\draw  [dash pattern={on 2.5pt off 2pt}] (2,2) -- (0,2);
			\draw  [dash pattern={on 2.5pt off 2pt}] (6,0) -- (6,2);
			\draw  [dash pattern={on 2.5pt off 2pt}] (2,0) -- (2,2);
            \draw  [dash pattern={on 2.5pt off 2pt}] (8,0) -- (8,4);
			\node at (2,-0.3) {\footnotesize $1$};
			\node at (6,-0.3) {\footnotesize $3$};
			\node at (8,-0.3) {\footnotesize $4$};
			\node at (0,-0.3) {\footnotesize $0$};
			\node at (-0.3,2) {\footnotesize $1/2$};
            \node at (-0.3,4) {\footnotesize $1$};
		\end{tikzpicture}
	\caption{The function $\eta_2(x)$}
			\label{Fig-250902}
	\end{minipage}
\end{figure}
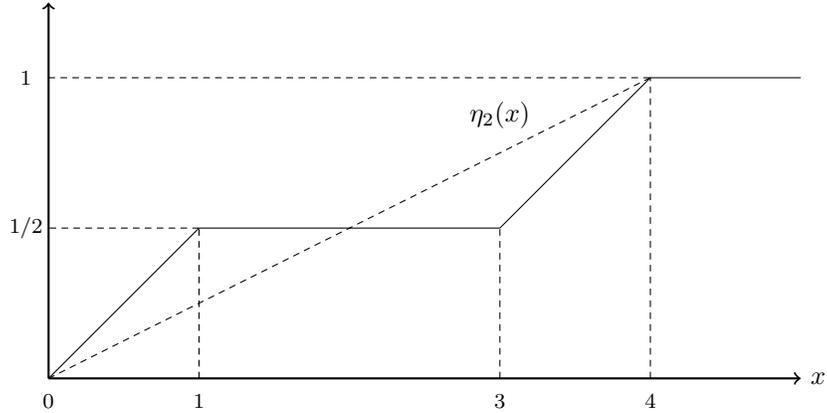

\begin{example}
 \label{ex-250921} {\rm
$(\mathcal{V}\varsubsetneq \mathcal{H}^\ast)$.\ \ Let $F_2$ be a distribution function with $F_2(0-)=0$, and $\eta_2(x)=\Fbar_2(1/x)$ be defined as follows (see Figure \ref{Fig-250902})
$$
   \eta_2(x) =\left\{ \begin{array}{ll} x/2, & x\in [0, 1],\\
       1/2, & x\in (1,3],\\
       x/2-1, & x\in (3, 4), \\
       1, & x\in [4, \oo).  \end{array}
   \right.
$$
It is easy to see that $\eta_2(x)/x=(1/x) \Fbar_2(1/x)$ is not decreasing in $x\in\R_{++}$, which implies $F_2\not\in\mathcal{V}$. Now, we prove that $\eta_2$ is subadditive on $\R_+$, that is,
\begin{equation}
  \label{eq-250902}
    \eta_2(x+y) \le \eta_2(x) +\eta_2(y),\quad x, y\in\R_{++}.
\end{equation}
Notice that
\begin{itemize}
 \item Since $\eta_2(z)/z$ is decreasing in $z\in (0,3]$, \eqref{eq-250902} holds true when $x+y\le 3$.
 \item When $x+y\in (3,\oo)]$ with $x\ge 1$ and $y\ge 1$, we have $\eta_2(x)+\eta_2(y)\ge 1/2+1/2\ge \eta_2(x+y)$. When $x+y\in (3,\oo)$ with $x\in (0,1]$, we have $\eta_2(x+y)-\eta_2(y)\le \eta_2(x)$.
\end{itemize}
Then \eqref{eq-250902} always holds, implying $F_2\in \mathcal{H}^\ast$. Therefore, $\mathcal{V} \varsubsetneq \mathcal{H}^\ast$.    }
\end{example}

Note that $F$ in Example~\ref{ex-250921} is a continuous distribution function. 
\cite{ALO24} provided a discrete distribution in their Example~2.7 that belongs to $\mathcal{H}^\ast$ but not to $\mathcal{H}$, because every distribution in $\mathcal{H}$ is continuous on $\mathbb{R}_{++}$ (see Remark~\ref{continue}).

\begin{example}[${\cal G} \not\subset {\cal H}$]
 \label{ex-250926}
{\rm Let $F$ be a Log-Cauchy distribution, that is,
\begin{equation*}
    F(x) = \frac{\arctan(\log x)}{\pi} + \frac{1}{2}, \quad x \in\R_{++}.
\end{equation*}
Then the density function of $F$ is
$$
   f(x) = \frac {1}{\pi x [1+(\log x)^2]}, \quad x\in\R_{++}.
$$
According to Table 1 of \cite{CHSZ25}, we have $F \in\mathcal{H}$. In Appendix \ref{Appendix-A}, it is shown that $F^\beta \notin \mathcal{H}$ for $\beta= 0.5$.

Next, we prove $F \in \mathcal{G}$, i.e., $\Lambda_F(x)=-\log F(1/x)$ is subadditive on $\R_{++}$. If so, by Proposition \ref{pr-250925} (ii), we have $F^\beta \in \mathcal{G}$ for all $\beta \in (0,1)$. To establish the subadditivity of $\Lambda_F$, it suffices to show that $L(x)= \Lambda_F(x)/x$ is decreasing on $\R_{++}$. In view of $F(1/x)=\Fbar (x)$ and $(1/x) f(1/x) = xf(x)$, we have
\begin{equation*}
    L'(x) = \frac{1}{x^2}\[\log \Fbar (x)) + \frac{x f(x)}{\Fbar (x)}\]\stackrel {\rm sgn}{=} \log \Fbar (x)) + \frac{x f(x)}{\Fbar (x)},
\end{equation*}
which is non-positive for all $x\in\R_{++}$ (For its proof, see Appendix \ref{Appendix-A}). Therefore, $F \in \mathcal{G}$.   }
\end{example}

Note that $F$ in Example~\ref{ex-250926} is a continuous distribution function. 
\cite{CS25} provided a discrete distribution in their Example~3 that belongs to $\mathcal{G}$ but not to $\mathcal{H}$, because every distribution in $\mathcal{H}$ is continuous on $\mathbb{R}_{++}$ (see Remark \ref{continue}).

\begin{example}[${\cal G} \not\subset {\cal V}$]
{\rm
Let $F$ be a distribution function such that $F(0-)= 0$ and
$$
    \Lambda_F(x) = -\log F\(\frac {1}{x}\)= \left\{\begin{array}{ll} x^{1/2}, & x \in [0, 1], \\
        (x-0.99)^{1/2} + 0.9, & x \ge 1.  \end{array} \right.
$$
We first show the subadditivity of $\Lambda_F$, i.e., $F \in \cal G$. Choose $x \ge 0$ and $y \ge 0$. If $x +y \le 1$, then $\Lambda_F(x+y) = (x+y)^{1/2} \le x^{1/2} + y^{1/2} = \Lambda_F(x) + \Lambda_F(y)$. If $x + y > 1$, we need to consider the following three cases.
\begin{itemize}
 \item Case 1. If $x \ge 1$ and $y \ge 1$, then $\Lambda_F(x+y) = (x+y -0.99)^{1/2} + 0.9 \le (x-0.99)^{1/2} + (y-0.99)^{1/2} + 1.8 = \Lambda_F(x) + \Lambda_F(y)$.

 \item Case 2. If $x\ge 1$ and $0\le y < 1$, then $\Lambda_F(x+y) = (x+y -0.99)^{1/2}+ 0.9 \le (x-0.99)^{1/2} + y^{1/2} + 0.9 = \Lambda_F(x) + \Lambda_F(y)$. The proof for the case $0 \le x < 1$ and $y \ge 1$ is similar.

 \item Case 3. If $0 \le x<1$ and $0 \le y< 1$, we have $x+y \in [1, 2)$ and
     $$
        \Lambda_F(x+y)= (x+y -0.99)^{1/2}+ 0.9 \le (x+y-1)^{1/2}+1
          \le x^{1/2} + y^{1/2} =\Lambda_F(x) + \Lambda_F(y).
     $$
\end{itemize}
Define $g(x)= (1/x) \Fbar (1/x)= [1-\exp\{-\Lambda_F(x)\}]/x$. It can be checked that $g(1) \approx 0.6321 < g(1.01) \approx 0.6406$. This means $g(x)$ is not decreasing. Thus, $F \notin \mathcal{V}$.   }
\end{example}

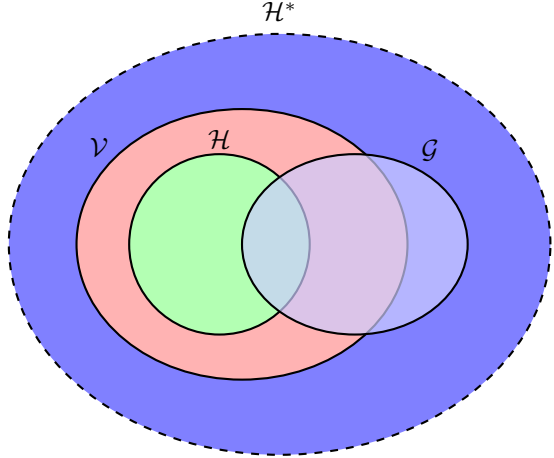
\begin{figure}[htbp]
\centering
\begin{tikzpicture}
\begin{scope}
    \draw[thick, fill=blue!50, dashed] (0,0) ellipse (3.6cm and 2.8cm);
    \node at (0,3.1) {$\mathcal{H}^*$};
\end{scope}
\begin{scope}
    \draw[thick, fill=red!30] (-0.5,0) ellipse (2.2cm and 1.8cm);
    \node at (-2.4,1.3) {$\mathcal{V}$};
\end{scope}
\begin{scope}
    \draw[thick, fill=green!30] (-0.8,0) circle (1.2cm);
    \node at (-0.8,1.4) {$\mathcal{H}$};
\end{scope}
\begin{scope}
    \draw[thick, fill=blue!20, fill opacity=0.7] (1,0) ellipse (1.5cm and 1.2cm);
    \node at (2,1.25) {$\mathcal{G}$};
\end{scope}
\end{tikzpicture}
    \caption{Venn diagram illustrating the relationships among the classes $\mathcal{H}$, $\mathcal{G}$, $\mathcal{V}$, and $\mathcal{H}^*$. The largest class $\mathcal{H}^*$ is indicated by the dashed boundary, while $\mathcal{H}$ is a subset of $\mathcal{V}$, and $\mathcal{G}$ has non-empty intersections with both $\mathcal{H}$ and $\mathcal{V}$.}
\label{Fig-Venn1}
\end{figure}

The above discussion thus allows us to depict the relationships in Figure \ref{Fig-Venn1} among the classes $\mathcal{H}$, $\mathcal{G}$, $\mathcal{V}$ and $\mathcal{H}^\ast$ in a Venn diagram.

Proposition \ref{pr-250925}(ii) shows that $F^\beta \in \mathcal{H}$ [resp. $\mathcal{V}$, $\mathcal{H}^\ast]$ when $F \in \mathcal{H}$ [resp. $\mathcal{V}$, $\mathcal{H}^\ast]$ and $\beta \ge 1$. Below, we demonstrate that this result cannot be extended to $\beta \in (0, 1)$.

\begin{example}
 \label{ex-250927}
{\rm Let $F$ be a distribution function with $F(0-) = 0$, and $\eta(x) = \Fbar(1/x)$ be defined as follows
\begin{equation*}
   \eta(x) = \left\{\begin{array}{ll} \dfrac{x}{2}, & x \in [0, 1], \\[8pt]
      \dfrac{1}{2} + \dfrac{x-1}{4}, & x \in (1, 3], \\[5pt] 1, & x \ge 3.  \end{array} \right.
\end{equation*}
Then $\eta(x)$ is a concave function on $\R_{++}$, i.e., $F \in \mathcal{H}$. Hence, $F \in \mathcal{V}$ and $F \in \mathcal{H}^\ast$. Now, define $G = F^\beta$ for $\beta \in (0, 1)$, so that $\eta_\beta(x) := \Gbar(1/x) = 1 - \(1 - \eta(x)\)^\beta.$ For $\beta \in (0, 0.69)$, we have $\eta_\beta(3) > \eta_\beta(2) + \eta_\beta(1)$, which implies that $F \notin {\cal H}^*$.  }
\end{example}

A counterexample is given in Example \ref{ex-250919} to show that the likelihood ratio order $\le_{\rm lr}$ in Proposition \ref{pr-250925} (iv) for $\mathcal{H}$ cannot be replaced by the hazard rate order $\le_{\rm hr}$. Specifically, there exist distributions $F$ and $G$ such that $F \in \cal H$ and $F \le_{\rm hr} G$, yet $G \notin \cal H$.

\begin{example}
  \label{ex-250919}
{\rm Consider two distribution functions $F$ and $G$, having a common support $(6/5, \oo)$, with survival functions given by
\begin{equation*}
   \Fbar(x) = \frac{c_1}{x+1}\ \ {\rm and}\ \ \Gbar(x)= c_2 \frac{3(1/x - 1)^2+1}{x}\ \ {\rm for}\ x \ge \frac {6}{5},
\end{equation*}
where the positive constants $c_1$ and $c_2$ are determined such that $\Fbar(6/5)= \Gbar(6/5)= 1$. To verify that $G$ is a distribution, it suffices to prove that $h(y):=\Gbar(1/y)=c_2 [3(y-1)^2+1] y$ is increasing in $y\in (0, 5/6)$. This is trivial since $h'(y)=c_2 (3y-2)^2\ge 0$.

It is easy to show that $\Fbar(1/x) = c_1 x / (1+x)$ is concave in $x\in\R_+$, and hence $F \in \cal H$. Note that
\begin{equation*}
   g(x) := \frac{\Gbar(1/x)}{\Fbar(1/x)} = \frac{c_2}{c_1} [3(x - 1)^2+1](1+x),\quad x \le \frac{5}{6}.
\end{equation*}
Since $g'(x) = 9x^2 - 6x - 2 = 9 (x - 1/3)^2 - 3 < 0$ for $x\in (0,5/6)$, we have $F \le_{\rm hr} G$. However, $\Gbar(1/x)$ is convex over $[2/3, 5/6]$. This means $G \notin \cal H$.    }
\end{example}

In Proposition \ref{pr-250925} (iv), $F\in\mathcal{G}$ and $F\le_{\rm hr} G$ does not imply $G\in\mathcal{G}$, as shown by the next example.

\begin{example}
{\rm Let
\begin{equation*}
    \Fbar(x) =\frac{1}{1+x}\ \ \hbox{and}\ \ \Gbar(x)=\min\left\{\frac{2}{1+x}, 1\right\}\ \hbox{for}\ x \in \R_+.
\end{equation*}
It is known that $F\in\mathcal{G}$ \citep[see][Example 2]{CS25}. Note that
\begin{equation*}
\frac{\Gbar(x)}{\Fbar(x)} = \left\{\begin{array}{ll} 1+x, & 0 \leq x < 1, \\ 2, & x \geq 1.  \end{array} \right.
\end{equation*}
Therefore, $F \leq_{\rm hr} G$. However, the subadditivity of $\Lambda_G$ does not hold in general, which can be checked by choosing $x= y=0.4$. This means $G \notin \mathcal{G}$. In fact, it is easy to see  $G \notin \mathcal{G}$ since $\essinf(G)=1$, not zero.  }
\end{example}

Many of the closure properties presented in the following propositions can be understood through the stochastic transform order characterizations given in Remark \ref{convex_trans}. For $\mathcal{H}$, the condition $U \leq_{\rm c} 1/X$ implies that applying a convex transformation to $X$ preserves the class, since convex transformations preserve the convex transform order. Analogously, $\mathcal{V}$ and $\mathcal{H}^\ast$ are closed under convex transformations due to the preservation of the anti-starshaped and subadditive orders, respectively. 
For completeness, we still provide direct proofs of these properties based on the original definitions of the classes.

\begin{proposition} \
\label{pr-250928}
\begin{itemize}
  \item[{\rm (vi)}] Let $X$ and $Y$ be independent.
    \begin{itemize}
    \item If $X, Y\in\mathcal{G}$, then $\max\{X, Y\}\in \mathcal{G}$ \citep[][Proposition 2]{CS25}.

    \item if $X, Y\in\mathcal{H}$, then $\max\{X, Y\}\in \mathcal{H}$ \citep[][Proposition 3]{CHSZ25}.

    \item If $X, Y\in\mathcal{V}$ [resp. $\mathcal{H}^\ast$], then $\max\{X, Y\}\in \mathcal{V}$ [resp. $\mathcal{H}^\ast$].
    \end{itemize}

  \item[{\rm (vii)}] If $X\in \mathcal{H}$, then $(X-c)_+\in\mathcal{H}$ for any $c\in\R_{++}$ \citep[][Proposition 5]{CHSZ25}. \\
       If $X\in \mathcal{V}$ [resp. $\mathcal{G}$, $\mathcal{H}^\ast$], then $(X-c)_+\in\mathcal{V}$ [resp. $\mathcal{G}$, $\mathcal{H}^\ast$] for any $c\in\R_{++}$.

  \item[{\rm (viii)}] Let $X$ and $Y$ be independent such that $X\in\mathcal{V}$ [resp. $\mathcal{H}$, $\mathcal{G}$, $\mathcal{H}^\ast$]. If $Y$ is non-negative, then $(X-Y)_+\in \mathcal{V}$ [resp. $\mathcal{H}$, $\mathcal{G}$, $\mathcal{H}^\ast$].

  \item[{\rm (ix)}]  If $X\in \mathcal{V}$ [resp. $\mathcal{H}$, $\mathcal{H}^\ast$], then $[X|X>c] \in \mathcal{V}$ [resp. $\mathcal{H}$, $\mathcal{H}^\ast$] for any $c\in\R_{++}$. However, $[X|X>c]\notin \mathcal{G}$ for any $c\in\R_{++}$.
\end{itemize}
\end{proposition}

The next examples demonstrate that $\mathcal{H}$, $\mathcal{V}$, $\mathcal{H}^\ast$ and $\mathcal{G}$ are not closed under convolution.

\begin{example}[Convolution]
{\rm Let $X_1$ and $X_2$ be iid Pareto$(1)$ distributed random variables. It is easy to see that $X_1\in \mathcal{H}$ and hence $X_1\in\mathcal{V}$ and $X_1\in\mathcal{H}^\ast$ by Proposition \ref{pr-250925}(i). We claim that $X_1 + X_2 \notin \mathcal{H}^\ast$ and hence $X_1+X_2\not\in \mathcal{H}$ and $X_1+X_2\not\in \mathcal{V}$. To see it, the distribution function of $X_1+X_2$ is given by
$$
    G(x)=\left \{ \begin{array}{ll} 0, & x\le 2,\\
        1 - 2 x^{-1} - 2 x^{-2} \log (x-1),  & x \ge 2. \end{array}\right.
$$
However, the inequality $\Gbar(1/(x+y)) \leq \Gbar(1/x) + \Gbar(1/y)$ does not hold in general for any $(x, y)\in\R_{++}$. A counterexample is given by $x =y= 0.1$. This means $X_1+X_2\notin \mathcal{H}^\ast$. Therefore, $\mathcal{H}$, $\mathcal{V}$ and $\mathcal{H}^\ast$ are not closed under convolution.
 }
\end{example}

\begin{example}[Convolution]
{\rm
Let $X, X_1, X_2$ be iid with distribution function
\begin{equation*}
   F(x) = \frac{x}{1+x},\quad x \in \R_+.
\end{equation*}
It is known $F\in\mathcal{G}$. However, $X_1 + X_2  \notin \cal G$. To prove it, denote $Z = X_1 + X_2 \sim G$. Then
\begin{align*}
  G(z) &= \P\(X_1 + X_2 \leq z\) = \int_0^z \int_0^{z-x} \frac{1}{(1+x)^2(1+y)^2} \d y \d x  \\
    &= \frac{z}{z+2} - \frac{2\log(1+z)}{(z+2)^2},
\end{align*}
and
$$
  \Lambda_G(x) = - \log G\(\frac{1}{x}\) = 2\log(1+2x) - \log\(1+2x-2x^2 \log\(1+ \frac{1}{x}\)\).
$$
Choosing $x = 0.02$ and $y = 0.18$, we have $\Lambda_G(x+y)=\Lambda_G(0.2) \approx 0.444488 > \Lambda_G(x)+ \Lambda_G(y) \approx 0.443596$. Thus, $G\notin\mathcal{G}$.    }
\end{example}

It is a common consensus that applying an increasing, convex and nonconstant transformation to a random variable $X$ results in a new random variable $Y$ with a heavier right tail than $X$. The following properties demonstrate that $\mathcal{H}$, $\mathcal{V}$, $\mathcal{H}^\ast$ and $\mathcal{G}$ are closed under an increasing, convex and nonconstant transform anchoring at zero:

\begin{itemize}
  \item[P${}_1$] {\rm \citep[][Lemma 5.5]{Vin25}.}\ \ Let $\psi$ be an increasing, convex and nonconstant function with $\psi(0)=0$. If $X\in \mathcal{V}$, then $\psi(X)\in \mathcal{V}$.

  \item[P${}_2$] {\rm \citep[][Proposition 3]{CHSZ25}.}\ \ Let $\psi$ be a strictly increasing with $\psi(0)=0$ and $1/\psi^{-1}(1/x)$ being concave in $x\in\R_{++}$. If $X\in \mathcal{H}$, then $\psi(X)\in \mathcal{H}$.

  \item[P${}_3$] {\rm \citep[][Theorem 2.9]{ALO24}.}\ \ Let $\psi$ be a continuous, and nonconstant star-shaped function with $\psi(0)=0$. If $X\in \mathcal{H}^\ast$, then $\psi(X)\in \mathcal{H}^\ast$.

  \item[P${}_4$] {\rm \citep[][Proposition 2(iv)]{CS25}.}\ \ Let $\psi$ be an increasing, convex and nonconstant function with $\psi(0)=0$. If $X\in \mathcal{G}$, then $\psi(X)\in \mathcal{G}$.
\end{itemize}

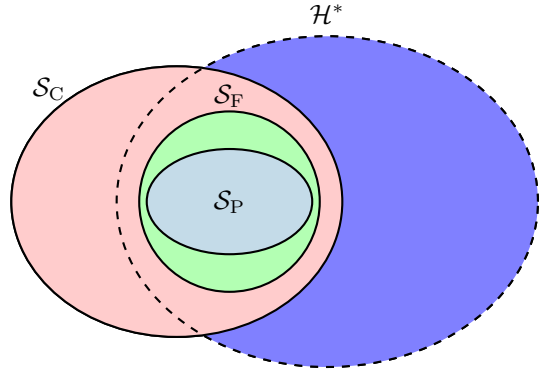
\begin{figure}[htbp]
\centering
\begin{tikzpicture}
\draw[thick, dashed] (0.5,0) ellipse (2.8cm and 2.2cm);
\node at (0.5,2.5) {$\mathcal{H}^*$};
\draw[thick] (-1.5,0) ellipse (2.2cm and 1.8cm);
\node at (-3.2,1.5) {$\mathcal{S}_{\mathrm{C}}$};
\begin{scope}
    \clip (0.5,0) ellipse (2.8cm and 2.2cm);
    \fill[blue!50] (-4,-3) rectangle (4,3);
\end{scope}
\begin{scope}
    \clip (-1.5,0) ellipse (2.2cm and 1.8cm);
    \fill[red!20] (-4,-3) rectangle (4,3);
\end{scope}
\draw[thick, dashed] (0.5,0) ellipse (2.8cm and 2.2cm);
\draw[thick] (-1.5,0) ellipse (2.2cm and 1.8cm);
\begin{scope}
    \draw[thick, fill=green!30] (-0.8,0) circle (1.2cm);
    \node at (-0.8,1.4) {$\mathcal{S}_{\mathrm{F}}$};
\end{scope}
\begin{scope}
    \draw[thick, fill=blue!20, fill opacity=0.7] (-0.8,0) ellipse (1.1cm and 0.7cm);
    \node at (-0.8,0) {$\mathcal{S}_{\mathrm{P}}$};
\end{scope}
\end{tikzpicture}
\caption{Venn diagram illustrating the relationships among four classes $\mathcal{S}_{\rm P}$, $\mathcal{S}_{\rm F}$, $\mathcal{S}_{\rm C}$ and $\mathcal{H}^\ast$.}
\label{Fig-Venn2}
\end{figure}

\begin{remark}
\label{re-20251031}
{\rm
Let $\psi$ be an increasing, convex and nonconstant function, and denote $Y=\psi(X)$. If $\psi(0)=0$ and $X$ has  Pareto(1) distribution, then we say $Y$ or its distribution is \emph{super-Pareto} \citep{CEW25}. If $\psi(0)=0$ and $X$ has Fr\'{e}chet(1) distribution, we say $Y$ or its distribution is \emph{super-Fr\'{e}chet} \citep{CS25}. If $\psi(-\oo)=0$ and $X$ has Cauchy$(0,1)$ distribution given by $F_{\rm C}(x)= \pi^{-1}\arctan (x) +1/2$ for $x\in\R$, we say $Y$ or its distribution is \emph{super-Cauchy} \citep{Mul25}. Denote by $\mathcal{S}_{\rm P}$, $\mathcal{S}_{\rm F}$ and $\mathcal{S}_{\rm C}$ the classes of all super-Pareto, super-Fr\'{e}chet and super-Cauchy distributions, respectively.

Denote by $\mathcal{F}_+$ the class of all distributions of non-negative random variables. Then
\begin{align*}
 & \mathcal{S}_{\rm P} =\{G\in \mathcal{F}_+:   \hbox{Pareto} (1) \le_{\rm c} G\}, \\
 & \mathcal{S}_{\rm F} =\{G\in \mathcal{F}_+:   \hbox{Fr\'{e}chet} (1) \le_{\rm c} G\}, \\
 & \mathcal{S}_{\rm C} =\{G\in \mathcal{F}_+:   \hbox{Cauchy} (0,1) \le_{\rm c} G\}.
\end{align*}
Since $\hbox{Fr\'{e}chet} (1) \le_{\rm c} \hbox{Pareto} (1)$, $\hbox{Pareto} (1) \not\le_{\rm c} \hbox{Fr\'{e}chet} (1)$ \citep[][Example 4]{CS25} and $\hbox{Cauchy} (0,1) \le_{\rm c} \hbox{Fr\'{e}chet} (1)$, $\hbox{Fr\'{e}chet} (1) \not\le_{\rm c} \hbox{Cauchy}(0,1)$ \citep[][Theorem 2.10]{Mul25}, we have
$$
    \mathcal{S}_{\rm P}\varsubsetneq \mathcal{S}_{\rm F}\varsubsetneq \mathcal{S}_{\rm C}.
$$
\cite{Mul25} gave counterexamples to show $\mathcal{G}\not\subset\mathcal{S}_{\rm C}$ and $\mathcal{S}_{\rm C}\not \subset \mathcal{H}^\ast$. It is easy to check that, for $\alpha\in (0,1]$,
$$
    \hbox{Pareto}(1) \le_{\rm c} \hbox{Pareto}(\alpha),\qquad \hbox{Fr\'{e}chet}(1) \le_{\rm c} \hbox{Fr\'{e}chet}(\alpha).
$$
Thus, $\hbox{Pareto}(\alpha)\in \mathcal{S}_{\rm P}$ and $\hbox{Fe\'{e}chet}(\alpha)\in \mathcal{S}_{\rm F}$ for $\alpha\in (0,1]$. By Property P${}_3$ and $\hbox{Fr\'{e}chet}(1)\in\mathcal{H}^\ast$, we have
$\mathcal{S}_{\rm F} \subset \mathcal{H}^\ast$. The relationships among four classes $\mathcal{S}_{\rm P}$, $\mathcal{S}_{\rm F}$, $\mathcal{S}_{\rm C}$ and $\mathcal{H}^\ast$ are depicted in Figure \ref{Fig-Venn2}.
  }
\end{remark}

Examples \ref{ex-250923} and \ref{ex-250924} below show that $\mathcal{H}$, $\mathcal{V}$ and $\mathcal{H}^\ast$ are not closed under a simple convex transform $\psi(x)=x+c$ with $c>0$. It is also shown that the assumption $\psi(0)=0$ cannot be removed from Properties P$_1$--P$_4$.

\begin{example}
  \label{ex-250923}
{\rm Let $X\sim F$, where $F$ is the Fr\'{e}chet$(1)$ distribution. Denote $Y=X+1=\psi(X)\sim G$, where $\psi(x)=x+1$. Let $f$ and $g$ denote the respective density functions of $X$ and $Y$. It is easy to see that $x^2f(x)$ is increasing on $\R_+$, while $x^2 g(x)$ is increasing on $[0,2]$ and decreasing on $(2, \oo)$. Therefore, $X\in\mathcal{H}$ while $Y\notin\mathcal{H}$.   }
\end{example}

\begin{example}
 \label{ex-250924}
{\rm Let $X\sim F_2$ with $F_2$ given by Example \ref{ex-250921}. Denote $Y=X+1=\psi(X)\sim G$, where $\psi(x)=x+1$. It is easy to see that
$$
  \eta_{{}_Y}(x):= \Gbar\(\frac {1}{x}\)  =\left\{ \begin{array}{ll} \dfrac {x}{2(1-x)}, & x\in \[0, \dfrac {1}{2}\],\\[8pt]
       \dfrac {1}{2}, & x\in \(\dfrac {1}{2}, \dfrac {2}{3}\],\\[8pt]
       \dfrac {x}{2(1-x)}-1, & x\in \(\dfrac {2}{3}, \dfrac {4}{5}\), \\[8pt]
       1, & x\in \[\dfrac {4}{5}, \oo\).  \end{array}
   \right.
$$
Choosing $x_0=y_0=2/5$, we have $\eta_{{}_Y}(x_0)+\eta_{{}_Y}(y_0)=2\eta_{{}_Y}(x_0)=2/3<1=\eta_{{}_Y}(x_0+y_0)$, violating the subadditivity of $\eta_{{}_Y}(x)$. Thus, $Y\notin\mathcal{H}^\ast$ and hence $Y\notin \mathcal{V}$. However, $X\in\mathcal{V}$ and hence $X\in\mathcal{H}^\ast$, as shown in Example \ref{ex-250920}.
}
\end{example}

The next example demonstrates that the class $\mathcal{H}$ [resp. $\mathcal{G}$, $\mathcal{V}$ and $\mathcal{H}^\ast$] is not closed under weak convergence.

\begin{example}
 \label{ex-250925}
{\rm Consider the distribution functions
\begin{equation*}
   F_n(x) = 1 - \frac{1}{nx+1},\quad x\in\R_+.
\end{equation*}
Note that $\Fbar_n(1/x) = x/(n+x)$ is concave on $\R_+$, so $F_n \in \mathcal{H}$ for all $n$. Hence, $F_n \in \mathcal{V}$ and $F_n \in \mathcal{H}^\ast$. Also, $\Lambda_F(x)=-\log F(1/x)$ is subadditive on $\R_+$, i.e., $F \in \mathcal{G}$. However, as $n \to \oo$, $F_n$ converges weakly to the degenerate distribution at zero, which does not belong to ${\cal H}^\ast$. }
\end{example}

\section{Stochastic dominance between diversified portfolios }
\label{sect-4}

\subsection{$\mathcal{H}$ and $\mathcal{H}^\ast$-distributed losses triggered by events}

In actuarial science, extremely heavy-tailed losses are often triggered by events with small probabilities of occurrence \citep{BGHJN97}. In this context, the outcome (loss) of a rare event can be modeled as $X \id_A$, where $X$ is a heavy-tailed random variable and $A$ is the triggering event independent of $X$. Let $\bm X=(X_1, \ldots, X_n)$ be a vector of $n$ iid random variables with a common distribution $F\in\mathcal{H}$, and $A_1, \ldots, A_n$ be the respective triggering events of $X_1, \ldots, X_n$ such that $A_1, \ldots, A_n$ are independent of $\bm X$.

If $A_1=\cdots=A_n$, then $\id_{A_1}, \ldots, \id_{A_n}$ shares a comonotonicity structure, a notion of the strongest positive dependence. In this special dependence structure, by Theorem \ref{th-H}, we have
\begin{equation}
 \label{eq-251001}
   \sum^n_{i=1} \eta_i \id_{A_i} X_i \le_{\rm st} \sum^n_{i=1} \th_i \id_{A_i} X_i
\end{equation}
for all $\bm \th, \bm\eta\in\R_+^n$ such that $\bm\th\preceq_{\rm m} \bm \eta$. Theorem \ref{th-251002} below shows that inequality \eqref{eq-251001} also holds for any events $A_1, \ldots, A_n$ with an arbitrary dependence structure and an equal probability of occurrence. \cite{CHWZ25} in their Theorem 2 established Theorem \ref{th-251002} for the case $F$ being a Pareto$(\alpha)$, where $\alpha\in (0,1]$ is the tail parameter.

\begin{theorem}
  \label{th-251002}
Let $\bm X=(X_1, \ldots, X_n)$ be a vector of $n$ iid random variables with a common distribution $F\in\mathcal{H}$, and $A_1, \ldots, A_n$ be events with equal probability, which are independent of $\bm X$. Then \eqref{eq-251001} holds for all $\bm \th, \bm\eta\in\R_+^n$ such that $\bm\th\preceq_{\rm m} \bm \eta$.
\end{theorem}

\begin{proof}
Assume that $\bm\th, \bm\eta\in \Delta_n$ and $\P(A_i)=p\in (0,1)$ for each $i$. Below, we first show \eqref{eq-251001} for the case $n=2$. For $\lm\in (0,1/2]$, define $S(\lm)=\P\(\lm \id_{A_1} X_1 + (1-\lm) \id_{A_2} X_2>x\)$ for $x\in \R_+$. It suffices to show that $S(\lm)$ is increasing in $\lm\in (0, 1/2]$. Note that
\begin{align*}
  S(\lm) & = \P(A_1A_2)\, \P\big (\lm X_1 + (1-\lm) X_2>x\big ) + \P(A_1A_2^c) \Fbar\(\frac {x}{\lm}\)
                  + \P(A_1^cA_2) \Fbar\(\frac {x}{1-\lm}\) \\
  & = \P(A_1A_2)\, \P\big (\lm X_1 + (1-\lm) X_2>x\big ) + (p-\P(A_1A_2)) \[\Fbar\(\frac {x}{\lm}\)
                  + \Fbar\(\frac {x}{1-\lm}\)\].
\end{align*}
Since $F\in\mathcal{H}$, we have $\Fbar(x/\lm) +\Fbar(x/(1-\lm))$ is increasing in $\lm\in (0,1/2]$ for $x\in\R_+$. On the other hand, by Theorem \ref{th-H}, $\P\big (\lm X_1 + (1-\lm) X_2>x\big )$ is also increasing in $\lm\in (0,1/2]$ for $x\in\R_+$. Thus, $S(\lm)$ is increasing in $\lm\in (0,1/2]$ for $x\in\R_+$. This proves \eqref{eq-251001} for $n=2$.

Next, we consider the case $n\ge 3$ and $\bm\th\prec_{\rm m} \bm\eta$ by using the same argument as that in the proof of Theorem 2 in \cite{CHWZ25} with a minor modification. By the nature of majorization \citep[see][Section 1.A.3]{MOA11}, there exist a finite number of vectors $\bm\th^{(0)}, \bm \th^{(1)}, \ldots, \bm \th^{(m)}$ in $\R_+^n$ such that $\bm\th =\bm \th^{(0)} \prec_{\rm m} \bm \th^{(1)}\prec_{\rm m}\cdots\prec_{\rm m} \bm\th^{(m)}=\bm\eta$, and for each $k\in [m]$, $\bm\th^{(k-1)}$ and $\bm\th^{(k)}$ differ only in two coordinates. Without loss of generality, assume that $\bm\th$ and $\bm\eta$ differ only in coordinates $k$ and $\ell$ with $k<\ell$. For $S\subseteq [n]$, let $B_S=\(\bigcap_{i\in S }A_i\)\cap \(\bigcap_{i\in S^c }A_i^c\)$. For $\bm\th\in \R_+^n$, we write
\begin{align}
 \label{eq-250913}
    \sum_{i=1}^n\theta_i \id_{A_i} X_i & = \sum_{S\subseteq[n]\backslash \{k,\ell\}}\id_{B_S}\sum_{i\in S}\theta_iX_i
             +\sum_{\{k,\ell\}\subseteq S\subseteq[n]}\id_{B_S}\sum_{i\in S}\theta_iX_i \nonumber \\
    & \quad +\sum_{\{k\}\subseteq S\subseteq[n]\backslash\{\ell\}}\id_{B_S}\sum_{i\in S}\theta_i X_i
        +\sum_{\{\ell\}\subseteq S\subseteq[n]\backslash\{k\}}\id_{B_S}\sum_{i\in S}\theta_iX_i.
\end{align}
It is clear that
\begin{equation}
 \label{eq-250914}
    \sum_{S\subseteq[n]\backslash\{k,\ell\}}\id_{B_S}\sum_{i\in  S} \th_i X_i = \sum_{S\subseteq[n]\backslash\{k,\ell\}}\id_{B_S}\sum_{i\in S}\eta_iX_i.
\end{equation}
By Theorem \ref{th-H}, we have
\begin{equation}
 \label{eq-250915}
    \sum_{\{k,\ell\}\subseteq S\subseteq[n]}\id_{B_S}\sum_{i\in S}\theta_iX_i\ge_{\rm st} \sum_{\{k,\ell\}\subseteq S\subseteq[n]}\id_{B_S}\sum_{i\in S}\eta_iX_i.
\end{equation}
Note that
\begin{align*}
     \sum_{\{k\}\subseteq S\subseteq[n]\backslash\{\ell\}}\id_{B_S}\sum_{i\in S}\th_i X_i
           =\sum_{ D\subseteq[n]\backslash\{k,\ell\}}\id_{A_k}\id_{A_\ell^c}\prod_{s\in D}\id_{A_s}\prod_{t\in ([n]\backslash \{k,\ell\})\backslash D}\id_{A_t^c}\(\th_k X_k+\sum_{i\in D}\th_i X_i\).
\end{align*}
Then
\begin{align}
 \label{eq-250916}
   & \sum_{\{k\}\subseteq S\subseteq[n]\backslash\{\ell\}}\id_{B_S}\sum_{i\in S}\theta_i X_i +\sum_{\{\ell\}\subseteq S\subseteq[n]\backslash\{k\}}\id_{B_S}\sum_{i\in S}\th_i X_i \nonumber \\
   & \qquad = \sum_{ D\subseteq[n]\backslash\{k,\ell\}}\prod_{s\in D}\id_{A_s}\prod_{t\in ([n]\backslash \{k,\ell\})\backslash D} \id_{A_t^c} \nonumber \\
   & \qquad \qquad \times \(\id_{A_k}\id_{A_\ell^c}\bigg (\th_k X_k+\sum_{i\in D} \th_i X_i\bigg )
      +\id_{A_k^c}\id_{A_\ell}\bigg (\th_\ell X_\ell+\sum_{i\in D}\th_i X_i\bigg ) \).
\end{align}
For $D\subseteq[n]\backslash\{k,\ell\}$, let $G$ denote the distribution function of $\sum_{i\in D} \th_i X_i$. For $s\in\R_+$, we have
\begin{align*}
    &\P\(\id_{A_k}\id_{A_\ell^c}\bigg (\th_kX_k+\sum_{i\in D}\th_i X_i\bigg )
            +\id_{A_k^c}\id_{A_\ell}\bigg (\th_\ell X_\ell+\sum_{i\in D}\th_i X_i\bigg )>s\)\\
      & \quad = \P(A_k\cap A_\ell^c)\(\P\bigg ( \th_k X_k+\sum_{i\in D}\th_i X_i>s\bigg )
             +\P\bigg (\th_\ell X_\ell+\sum_{i\in D}\th_i X_i>s\bigg ) \)\\
      &\quad = \P(A_k\cap A_\ell^c)\(\int_{-\oo}^\oo \P(\th_k X_k>s-t)\,\d G(t)
              +\int_{-\oo}^\oo \P(\th_\ell X_\ell>s-t)\,\d G(t)\)\\
      &\quad = \P(A_k\cap A_\ell^c)\int_{-\oo}^\oo \[\Fbar\(\frac {s-t}{\th_k}\)
                +\Fbar\(\frac {t-s}{\th_\ell}\)\] \,\d G(t)\\
     &\quad \ge \P(A_k\cap A_\ell^c)\int_{-\oo}^\oo \[\Fbar\(\frac {s-t}{\eta_k}\)
                +\Fbar\(\frac {t-s}{\eta_\ell}\)\] \,\d G(t)\\
     & \quad= \P\(\id_{A_k}\id_{A_\ell^c}\bigg (\eta_k X_k+\sum_{i\in D}\eta_i X_i\bigg )
            +\id_{A_k^c}\id_{A_\ell}\bigg (\eta_\ell X_\ell+\sum_{i\in D}\eta_i X_i\bigg )>s\),
\end{align*}
where the inequality follows since $F\in\mathcal{H}$ and $(\th_k, \th_\ell) \prec_{\rm m} (\eta_k, \eta_\ell)$. From \eqref{eq-250916}, it follows that
\begin{align}
 \label{eq-250917}
   & \sum_{\{k\}\subseteq S\subseteq[n]\backslash\{\ell\}}\id_{B_S}\sum_{i\in S}\theta_iX_i
      +\sum_{\{\ell\}\subseteq S\subseteq[n]\backslash\{k\}}\id_{B_S}\sum_{i\in S}\theta_iX_i \nonumber \\
   & \qquad \ge_{\rm st} \sum_{\{k\}\subseteq S\subseteq[n]\backslash\{\ell\}}\id_{B_S}\sum_{i\in S}\eta_i X_i+\sum_{\{\ell\}\subseteq S\subseteq[n]\backslash\{k\}}\id_{B_S}\sum_{i\in S}\eta_i X_i.
\end{align}
Combining \eqref{eq-250913}-\eqref{eq-250915} and \eqref{eq-250917}, we conclude \eqref{eq-251001} for $n\ge 3$. This completes the proof of the theorem.
\end{proof}

Similarly, we can establish the next result.

\begin{theorem}
  \label{th-251003}
Let $\bm X=(X_1, \ldots, X_n)$ be a vector of $n$ iid random variables with a common distribution $F$, and $A_1, \ldots, A_n$ be events with equal probability, which are independent of $\bm X$. If $F\in\mathcal{H}^\ast$, then, for all $\bm \th\in \Delta_n$,
\begin{equation*}
    \id_{A_1} X_1 \le_{\rm st} \sum^n_{i=1} \th_i \id_{A_i} X_i.
\end{equation*}
\end{theorem}

\subsection{Losses with $\mathcal{H}$-type and $\mathcal{H}^\ast$-type tails}

In practice, random variables may not follow distributions from $\mathcal{H}$ or $\mathcal{H}^\ast$ in their entire support, whereas they have $\mathcal{H}$ or $\mathcal{H}^\ast$-type distributions beyond some thresholds.
Let $Y$ be a random variable with distribution function $G$ and $G(0-)=0$. We say that $Y$ has a $\mathcal{H}$-type distribution in tail beyond a point $c\in\R_{++}$ if there exists $F\in\mathcal{H}$ such that $\Gbar(y) =\Fbar(y)$ for $y\ge c$. Similarly, we can define a $\mathcal{H}^\ast$-type distribution in tail beyond a point $c$.

\begin{proposition}
  \label{pr-251002}
Let $Y_1, Y_2$ be iid random variables with distribution function $G$.

\begin{itemize}
 \item[{\rm (i)}] If $G$ is a $\mathcal{H}$-type distribution in tail beyond a point $c\in\R_{++}$, then
   $$
        \P\(\th_1 Y_1 +\th_2 Y_2 >x \) \ge \P\(\eta_1 Y_1+\eta_2 Y_2 > x \),\quad x\ge c,
   $$
   for $\bm \th, \bm\eta\in\Delta_2$ such that $\bm\th\prec_{\rm m} \bm\eta$.

 \item[{\rm (ii)}] If $G$ is a $\mathcal{H}^\ast$-type distribution in tail beyond a point $c\in\R_{++}$, then for $\bm \th\in\Delta_2$,
   $$
        \P\(\th_1 Y_1 +\th_2 Y_2 >x \) \ge \P\(Y_1>x \),\quad x\ge c.
   $$
\end{itemize}
\end{proposition}

\begin{proof}
We give the proof of part (i); the proof of part (ii) is similar. Assume that there exists $F\in\mathcal{H}$ such that $\Gbar(y)=\Fbar(y)$ for $y\ge c$, and let $X_1, X_2$ be iid with distribution function $F$. For $\lm\in (0,1/2]$, define $S(\lm)=\P(\lm Y_1+ (1-\lm) Y_2>x)$, where $x\ge c$. It suffices to show that $S(\lm)$ is increasing in $\lm\in (0, 1/2]$. Note that
\begin{align*}
  S(\lm) & = \P(\lm Y_1+ (1-\lm) Y_2>x, Y_1\le c) + \P(\lm Y_1+ (1-\lm) Y_2>x, Y_2\le c) \\
  & \qquad + \P(\lm Y_1+ (1-\lm) Y_2>x, Y_1> c, Y_2>c) \\[3pt]
  & =\int^c_0 \Gbar\(\frac {x-y}{\lm} +y\) \d G(y) +\int^c_0 \Gbar\(\frac {x-y}{1-\lm} +y\) \d G(y) \\[3pt]
  & \qquad + \P(\lm Y_1+ (1-\lm) Y_2>x, Y_1> c, Y_2>c) \\[3pt]
  & \stackrel {\rm def}{=} S_1(\lm) + S_2(\lm),
\end{align*}
where
\begin{align*}
  S_1(\lm) & =\int^c_0\[\Gbar\(\frac {x-y}{\lm}+y\) + \Gbar\(\frac {x-y}{1-\lm}+y\)\] \d G(y),\\[3pt]
  S_2(\lm) & =\P(\lm Y_1+ (1-\lm) Y_2>x, Y_1> c, Y_2>c).
\end{align*}
Since $(x-y)/\lm +y \ge c$ and $(x-y)/(1-\lm)+y\ge c$ for $y\in [0,c]$, we have
$$
    S_1(\lm) =\int^c_0\[\Fbar\(\frac {x-y}{\lm}+y\)+ \Fbar\(\frac {x-y}{1-\lm}+y\)\] \d G(y).
$$
It is shown in the proof of Theorem 1 in \cite{CHSZ25} that
$$
   \Fbar\(\frac {x-y}{\lm_2}+y\)+ \Fbar\(\frac {x-y}{1-\lm_2}+y\)
    \ge \Fbar\(\frac {x-y}{\lm_1}+y\)+ \Fbar\(\frac {x-y}{1-\lm_1}+y\)
$$
for $0<\lm_1<\lm_2\le 1/2$. Hence, $S_1(\lm)$ is increasing in $\lm\in (0,1/2]$.

Let $X_1^\ast, X_2^\ast$ be iid random variables with $X_1^\ast \stackrel {d}{=} [X_1|X_1>c]$. By Proposition
\ref{pr-250928}(ix), $X_1\in \mathcal{H}$ implies $X_1^\ast\in \mathcal{H}$. Then, by Theorem \ref{th-H},
\begin{align*}
  S_2(\lm) & = \P(\lm Y_1+ (1-\lm) Y_2>x | Y_1> c, Y_2>c) \[\Fbar(c)\]^2 \\
    & = \P(\lm X_1^\ast + (1-\lm) X_2^\ast >x ) \[\Fbar(c)\]^2,
\end{align*}
which is increasing in $\lm\in (0,1/2]$. Therefore, $S(\lm)=S_1(\lm) + S_2(\lm)$ is increasing in $\lm\in (0,1/2]$. This completes the proof of the proposition.
\end{proof}

Proposition \ref{pr-251002} may not hold for comparing the survival functions of $\sum^n_{i=1}\th_i Y_i$ and $\sum^n_{i=1} \eta_i Y_i$ for $n>2$. There is one gap in the proof of Proposition 3 in \cite{CHWZ25} for $n>2$, in which they considered Pareto-type distribution in tail beyond a point $c$.

\subsection{Truncated $\mathcal{H}$-distributed random variables}

Heavy-tailed distributions are widely used in finance and insurance due to their ability to capture extreme events. However, the infinite upper bound may raise concerns about theoretical practicality, as real word  risks often have natural limits. Truncated heavy-tailed distributions offer a more realistic approach by imposing an upper bound while retaining tail risk characteristics. For a threshold $c\in\R_{++}$, let $\bm Y=(X_1\wedge c, \ldots, X_n\wedge c)$ be a vector of the truncated random variables of $X_1, \ldots, X_n$ at $c$, where $X_1, \ldots, X_n$ are iid with a common distribution in $\mathcal{H}$. As the $Y_i$ have finite mean, one cannot expect to establish the usual stochastic ordering between $\sum^n_{i=1}\th_i Y_i$ and $\sum^n_{i=1} \eta_i Y_i$ for any $\bm\th,\bm\eta\in \Delta_n$ such that $\bm\th\prec_{\rm m} \bm\eta$. However, a more diversified portfolio $\sum^n_{i=1}\th_i Y_i$ can dominate a less diversified one $\sum^n_{i=1} \eta_i Y_i$ in the sense of tail probability in a large region if the upper bound $c$ is large enough.

\begin{proposition}
  \label{pr-250930}
Let $\bm\th, \bm\eta\in\Delta_n$ such that $\bm\th \prec \bm \eta$, and denote $b=1/\eta_{(1)}$, where $\eta_{(1)}=\min\{\eta_1, \ldots, \eta_n\}>0$. Let $X_1, \ldots, X_n$ be iid with a common distribution function $F\in \mathcal{H}$, and define $\bm Y=(X_1\wedge c, \ldots, X_n\wedge c)$ with $c\in (b, \oo)$. Then
$$
   \P\(\sum^n_{i=1} \eta_i Y_i >x\) \le \P\(\sum^n_{i=1} \th_i Y_i > x\),\quad x\in \[0, \frac {c}{b}\).
$$
\end{proposition}

\begin{proof}
The proof is similar to that of Proposition 6 in \cite{CHWZ25}. First, note that if there exists at least one $X_j>c$ with $j\in [n]$, then $\sum^n_{i=1} \eta_i (X_i\wedge c)\ge \eta_j X_j\ge \eta_{(1)} c=c/b$. Thus, for $x\in [0, c/b)$, we have
\begin{align*}
  \P\(\sum^n_{i=1} \eta_i Y_i \le x\) & = \P\(\sum^n_{i=1} \eta_i Y_i \le x,\ X_1\le c, \ldots, X_n\le c\) \\
  & = \P\(\sum^n_{i=1} \eta_i X_i \le x,\ X_1\le c, \ldots, X_n\le c\) \\
  & = \P\(\sum^n_{i=1} \eta_i X_i \le x\).
\end{align*}
Since $\bm\th\prec_{\rm m} \bm\th$, we have $c\eta_{(1)}\le c \th_{(1)}$. Similarly, for $x\in [0, c/b)$, we have
$$
   \P\(\sum^n_{i=1} \th_i Y_i \le x\) = \P\(\sum^n_{i=1} \th_i X_i \le x\).
$$
Hence, the desired result follows from Theorem \ref{th-H}.
\end{proof}

\subsection{$\mathcal{V}$-distributed losses}\label{subsec44}

In Theorem \ref{th-H}, ({\red SD}) is established under the assumption $F\in\mathcal{H}$. Since $\mathcal{H}\subset \mathcal{V}$, it is natural to wonder whether ({\red SD}) is also true if $F\in\mathcal{V}$. However, this assertion is negative, as shown by the following example.

\begin{example}
 \label{ex-250922}
{\rm Consider a random variable $X$ with survival function
\begin{equation*}
   \Fbar(x) = \left\{\begin{array}{ll} 1, & {\rm for}\ x < 1, \\
       1/x, & {\rm for}\ 1 \le x < 2, \\
       1/2, & {\rm for}\ 2 \le x < 3, \\
       3/(2x), & {\rm for}\ x \ge 3.  \end{array} \right.
\end{equation*}
It is easy to see that $x \Fbar(x)$ increases in $x\in\R_+$ and thus $F \in \cal V$. The corresponding density function is
\begin{equation*}
    f(x) = \left\{\begin{array}{ll} 1/x^2, & {\rm for}\ 1 \le x < 2, \\
         3/(2x^2), & {\rm for}\ x \ge 3,\\
         0, & {\rm otherwise}.  \end{array} \right.
\end{equation*}
Denote $A= [1, 2)$ and $B= [3, \oo)$, and let $X, X_1, X_2$ be iid random variables. Then,
\begin{align*}
   \P\(\frac{1}{4} X_1 + \frac{3}{4}X_2 > \frac{3}{2}\) &= \P\(X_1 + 3X_2 > 6\)  \\
      &= \P\(X_1\!+\! 3X_2 > 6, X_1\in\! A, X_2\in\! A\) + \P\(X_1\!+\! 3X_2 > 6, X_1 \in A, X_2 \in B\) \\
      & \quad +\P\(X_1\!+\!3X_2>6, X_1\!\in\! B, X_2\!\in\! A\) + \P\(X_1\!+\!3X_2>6, X_1\!\in\! B, X_2\!\in\! B\) \\
      &= \P\(X_1+ 3X_2 > 6, X_1 \in A, X_2 \in A\) + \P\(X_1 \in A, X_2 \in B\) \\
        &\quad + \P\(X_1 \in B, X_2 \in A\) + \P\(X_1 \in B, X_2 \in B\) \\
      &= 1 - \P\(X_1 + 3X_2 \leq 6, X_1 \in A, X_2 \in A\)  \\[2pt]
      &= 1 - \int_1^2 \frac {1}{x_1^2} \( \int_1^{(6-x_1)/3} \frac{1}{x_2^2} \d x_2 \) \d x_1  \\[2pt]
      &= 1 - \[\frac{1}{4} - \frac{1}{12} \ln\(\frac{5}{2}\)\] = \frac{3}{4} + \frac{1}{12} \ln\(\frac{5}{2}\) \approx 0.826358.
\end{align*}
Similarly,
\begin{align*}
  \P\(\frac{2}{5}X_1 + \frac{3}{5}X_2 > \frac{3}{2}\) &= \P\(2X_1 + 3X_2 > \frac {15}{2}\)  \\
     &= 1-\P\(2X_1 + 3X_2 \leq 7.5, X_1 \in A, X_2 \in A\)  \\[2pt]
     &= 1-\int_1^2\frac {1}{x_1^2} \(\int_1^{(7.5-2x_1)/3} \frac{1}{x_2^2} \d x_2\) \d x_1  \\[2pt]
     &= 1- \[\frac{3}{10} - \frac{8}{75} \ln\(\frac{22}{7}\)\]
        =\frac{7}{10}+\frac{8}{75}\ln\(\frac{22}{7}\) \approx 0.822147.
\end{align*}
It is known that $(2/5, 3/5) \preceq (1/4, 3/4)$. However, we observe that
\begin{equation*}
\P\(\frac{1}{4}X_1 + \frac{3}{4}X_2 > \frac{3}{2}\) > \P\(\frac{2}{5}X_1 + \frac{3}{5}X_2 > \frac{3}{2}\),
\end{equation*}
which implies
\begin{equation*}
    \frac{1}{4}X_1 + \frac{3}{4}X_2 \not \le_{\rm st} \frac{2}{5}X_1 + \frac{3}{5}X_2.
\end{equation*}   }
\end{example}

Theorem \ref{th-V} is a special consequence of one-basket-theorem in \cite{Vin25}, which 
was proved by applying the law of total probability and exploiting the special partition structure of the sample space $\Omega$. In the remaining of this subsection, we present a simple proof by the induction method.
\medskip

\noindent \emph{Proof of Theorem \ref{th-V}}.\ \ First, we consider the case $n=2$. For $(\th_1, \th_2)\in \Delta_2$ and $x\in\R_+$, we have
\begin{align*}
  \P(\th_1 X_1+\th_2 X_2>x) & \ge \P\(\th_1 X_1+\th_2 X_2>x, X_1>\frac {x}{\th_1}, X_2\le x\)   \\
    & \qquad + \P\(\th_1 X_1+\th_2 X_2>x, X_2>\frac {x}{\th_2}, X_1\le x\)   \\
    & \qquad + \P\(\th_1 X_1+\th_2 X_2>x, X_1>x, X_2> x\)  \\
   & = \P\(X_1>\frac {x}{\th_1}, X_2\le x\) + \P\(X_2> \frac {x}{\th_2}, X_1\le x\)
    + \P\(X_1>x, X_2>x\)  \\
   & =\Fbar_1\(\frac {x}{\th_1}\) F_2(x)+\Fbar_2\(\frac {x}{\th_2}\) F_1(x)+\Fbar_1(x)\Fbar_2(x)  \\
   & \ge \th_1\Fbar_1(x) F_2(x)+ \th_2\Fbar_2(x) F_1(x)+\Fbar_1(x)\Fbar_2(x)   \\
   & = \th_1 \Fbar_1(x) +\th_2\Fbar_2(x)\ge 0,
\end{align*}
where the last inequality follows from $F_1, F_2\in\mathcal{V}$. Now, assume $({\rm SD}_{\rm cp})$ holds when $n=m\ge 2$. For $(\th_1, \ldots, \th_m, \th_{m+1})\in\Delta_{m+1}$ and $x\in\R_+$, we have
\begin{align*}
  \P\(\sum^{m+1}_{i=1} \th_i X_i>x\) & =\P\( X_{m+1} > \frac {x}{\th_{m+1}}\)
        + \P\(\sum^{m+1}_{i=1} \th_i X_i>x, X_{m+1} \le \frac {x}{\th_{m+1}}\) \\[3pt]
  & = \Fbar_{m+1}\(\frac {x}{\th_{m+1}}\) + \int^{x/\th_{m+1}}_0 \P\(\sum^m_{i=1} \th_i X_i>x-\th_{m+1} t\)\d F_{m+1}(t) \\[3pt]
  & = \Fbar_{m+1}\(\frac {x}{\th_{m+1}}\) +\int^{x/\th_{m+1}}_0 \P\(\sum^m_{i=1} \frac {\th_i}{1-\th_{m+1}} X_i> \frac {x-\th_{m+1} t}{1-\th_{m+1}} \)\d F_{m+1}(t) \\[3pt]
  & \ge \Fbar_{m+1}\(\frac {x}{\th_{m+1}}\) +\int^{x/\th_{m+1}}_0 \sum^m_{i=1} \frac {\th_i}{1-\th_{m+1}} \Fbar_i\(\frac {x-\th_{m+1} t}{1-\th_{m+1}}\) \d F_{m+1}(t) \\[3pt]
  & = \sum^m_{i=1} \frac {\th_i}{1-\th_{m+1}}\[ \Fbar_{m+1}\(\frac {x}{\th_{m+1}}\) +\int^{x/\th_{m+1}}_0 \Fbar_i\(\frac {x-\th_{m+1} t}{1\!-\!\th_{m+1}}\) \d F_{m+1}(t)\] \\[3pt]
  & =\sum^m_{i=1} \frac {\th_i}{1\!-\!\th_{m+1}}\[ \P\(\frac {X_{m+1}}{\th_{m+1}} >x\) + \P\(\! (1\!-\!\th_{m+1}) X_i+\th_{m+1} X_{m+1} > x, \frac {X_{m+1}}{\th_{m+1}} \le x\) \] \\[3pt]
  & = \sum^m_{i=1} \frac {\th_i}{1-\th_{m+1}} \P\big ((1\!-\!\th_{m+1}) X_i+\th_{m+1} X_{m+1} > x\big ) \\[3pt]
  & \ge \sum^m_{i=1}\frac {\th_i}{1-\th_{m+1}}\big [(1-\th_{m+1})\Fbar_i(x)+\th_{m+1}\Fbar_{m+1}(x)\big ] \\[3pt]
  & =\sum^{m+1}_{i=1} \th_i \Fbar_i(x),
\end{align*}
where the first inequality follows from the induction assumption since $(\th_1/(1-\th_{m+1}), \ldots, \th_m/(1-\th_{m+1}))\in \Delta_m$, and the last inequality follows from the result for $n=2$.
This means $({\rm SD}_{\rm cp})$ holds when $n=m+1$. Therefore, the desired result follows by induction.

\section{Compound distributions}
\label{sect-5}

Let $\{Z_1,Z_2,\ldots\}$ be a sequence of iid random variables with distribution $F$, $N$ follow a Poisson distribution with parameter $\lambda\in \R_{++}$, and $N$ is independent of the $Z_i$. Then we say that $Y=\sum_{i=1}^N Z_i$ follows a \emph{compound Poisson distribution} with Poisson parameter $\lambda$ and distribution $F$, denoted by $C_{\rm Poi}(\lambda,F)$. Similarly, if $N\sim {\rm B}(m,p)$ [resp. ${\rm NB}(\alpha, p)$], then the distribution of $Y$ is called compound binomial distribution [resp. compound negative binomial distribution], denoted by $C_{\rm b}(m,p; F)$ [resp. $C_{\rm nb}(\alpha,p;F)$], where $m\in\N$, $\alpha\in\R_{++}$ and $p\in (0,1)$.

If $X_1, \ldots, X_n$ be iid $\sim F$, satisfying ({\red SD}) or ({\red ${\rm SD}^*$}), we also say $F$ satisfies ({\red SD}) or ({\red ${\rm SD}^*$}). It is known from \cite{CHSZ25} that
\begin{itemize}
  \item $C_{\rm Poi}(\lm, F)$ satisfies ({\red SD}) for any $\lm\in\R_{++}$ if and only if $F\in\mathcal{H}$;
    \item $C_{\rm Poi}(\lm, F)$ satisfies ${\rm ({\red {\rm SD}^\ast})}$ for any $\lm\in\R_{++}$ if and only if $F\in\mathcal{H}^\ast$.
\end{itemize}

\begin{theorem}\
 \label{th-251001}
Let $m\in \N$ be fixed with $m \ge 2$.
\begin{itemize}
 \item[{\rm (i)}] $C_{\rm b}(m,p; F)$ satisfies {\rm ({\red SD})} for any $p\in (0,1)$ if and only id $F\in \mathcal{H}$.
 \item[{\rm (ii)}] $C_{\rm b}(n,p; F)$ satisfies {\rm ({\red ${\rm SD}^*$})} for any $p\in (0,1)$ if and only if $F\in \mathcal{H}^*$.
\end{itemize}
\end{theorem}

\begin{proof}
We give the proof of part (i) by applying Theorem \ref{th-251002}; the proof of part (ii) is similar by applying Theorem \ref{th-251003}.

\emph{Sufficiency}\ \ Assume $F\in \mathcal{H}$. Using the argument similar to the proof of Theorem \ref{th-251002}, it suffices to establish ({\red SD}) for $n=2$. Let $Y_1, Y_2$ be iid random variables, each having $C_{\rm b}(m,p; F)$ distribution. If $F \in \mathcal{H}$, we need to show
\begin{equation*}
   \eta_1 Y_1 +\eta_2 Y_2 \le_{\rm st} \th_1 Y_1 + \th_2 Y_2
\end{equation*}
for $\bm\th, \bm\eta\in \Delta_2$ such that $\bm \th\prec_{\rm m} \bm\eta$.

First, we give a stochastic representation of a random variable $Y\sim C_{\rm b}(m,p; F)$. Denote by $\psi_Z(t)$ and $\psi_Y(t)$ the characteristic functions of $Z\sim F$ and $Y$, respectively. Note that $Y=\sum^N_{k=1} Z_k$, where $Z_1, \ldots, Z_m$ are iid with $Z_1\sim F$, and $N\sim {\rm B}(m,p)$, which is independent of the $Z_i$. Then the characteristic function of $Y$ is given by
\begin{align*}
   \psi_Y(t)=\E\[\exp\left\{ {\rm i}\, t \sum^N_{k=1} Z_k\right\} \] =\sum^m_{k=0} \[\psi_Z(t)\]^k {m\choose k} p^k (1-p)^{m-k} = \[1-p +p \psi_Z(t)\]^m,\quad t\in\R,
\end{align*}
which implies
\begin{equation}
\label{eq-250918}
      Y\ \stackrel {d}{=}\ \sum^m_{k=1} I_k Z_k,
\end{equation}
where $I_1, \ldots, I_m$ are iid ${\rm B}(1,p)$-distributed random variables, independent of the $Z_i$

Next, let $\big\{X_k^{(1)}, X_k^{(2)}, k\in [m]\big \}$ be iid random variables with a common distribution $F\in \mathcal{H}$, and let $\big \{I_k^{(1)}, I_k^{(2)}, k\in [m]\big \}$ be iid ${\rm B}(1,p)$-distributed random variables, independent of the $X_k^{(1)}$ and $X_k^{(2)}$. In view of \eqref{eq-250918}, we have
\begin{equation}
 \label{eq-250919}
  (Y_1, Y_2) \ \stackrel {d}{=}\ \(\sum^m_{k=1} I_k^{(1)} X_k^{(1)},\ \sum^m_{k=1} I_k^{(2)} X_k^{(2)}\).
\end{equation}
Thus,
\begin{align*}
  \eta_1 Y_1+\eta_2 Y_2\ & \stackrel {d}{=}\ \sum^m_{k=1} \(\eta_1 I_k^{(1)} X_k^{(1)}+\eta_2 I_k^{(2)} X_k^{(2)}\) \\
  & \le_{\rm st} \sum^m_{k=1} \(\th_1 I_k^{(1)} X_k^{(1)}+\th_2 I_k^{(2)} X_k^{(2)} \)\
   \stackrel {d}{=} \ \th_1 Y_1+\th_2 Y_2,
\end{align*}
where the inequality follows from Theorem \ref{th-251002} and the independence of all random variables. This proves part (i).

\emph{Necessity}\ \ In view of \eqref{eq-250919} and the independence of all random variables, we have, for any $x\in\R_+$.
\begin{align}
 \label{eq-250920}
  \P(\eta_1 Y_1+\eta_2 Y_2>x) & =\P\(\sum^m_{k=1}\eta_1 I_k^{(1)} X_k^{(1)}+\eta_2 I_k^{(2)} X_k^{(2)} >x\)\nonumber  \\
  & = m (1-p)^{2m-1} p \[\Fbar\(\frac {x}{\eta_1}\) + \Fbar \(\frac {x}{\eta_2}\)\] +\circ (p)\nonumber \\
  & = m p \[\Fbar\(\frac {x}{\eta_1}\) + \Fbar \(\frac {x}{\eta_2}\)\] +\circ (p),\quad p\to 0.
\end{align}
Similarly,
\begin{align}
  \label{eq-250921}
  \P(\th_1 Y_1+\th_2 Y_2>x) & =mp\[\Fbar\(\frac {x}{\th_1}\)+\Fbar \(\frac {x}{\th_2}\)\]+\circ (p),\quad p\to 0.
\end{align}
For any $\bm \th, \bm\eta\in \Delta_2$ satisfying $(\th_1, \th_2)\prec_{\rm m} (\eta_1, \eta_2)$, inequality ({\red SD}) implies $\P(\eta_1 Y_1+\eta_2 Y_2>x) \le \P(\th_1 Y_1+\th_2 Y_2>x)$ for $x\in\R_+$. Hence, letting $p\to 0$ in \eqref{eq-250920} and \eqref{eq-250921}, we have
$$
  \Fbar\(\frac {x}{\eta_1}\) + \Fbar \(\frac {x}{\eta_2}\) \le \Fbar\(\frac {x}{\th_1}\)+\Fbar \(\frac {x}{\th_2}\),\quad x\in\R_+.
$$
This means $F\in\mathcal{H}$.
\end{proof}

It is still unknown whether Theorem \ref{th-251001} holds for $C_{\rm nb}(\alpha, p;F)$.

\begin{remark} {\rm
If $C_{\rm b}(m,p; F)\in\mathcal{V}$ for any $p\in (0,1)$, then $F\in \mathcal{V}$. To see it, denote by $G_p(x)$ the distribution function of $C_{\rm b}(m,p; F)$, Then, for any $x\in\R_+$,
$$
  \Gbar_p(x)=\sum^m_{k=1} {m\choose k} p^k (1-p)^{m-k} \overline{F^{\ast k}}(x) = m p \Fbar(x) + \circ(p).
$$
Thus,
$$
   x \Gbar_p (x) = m p \cdot x\Fbar(x) + \circ(p), \quad p\to 0.
$$
So, if $x\Gbar_p(x)$ is increasing in $x\in\R_+$, we have $x\Fbar(x)$ is also increasing in $x\in\R_+$, i.e., $F\in\mathcal{V}$.  }
\end{remark}

\section{Discussions}
\label{sect-6}

For a random variable $X$ with distribution $F_X$, the VaR (Value-at-Risk) of $X$ at confidence level $\alpha\in [0,1]$ is defined to be the left inverse of its distribution function $F_X$, given by
$\VaR_\alpha(X):= F_X^{-1}(\alpha)$. We say that VaR is subadditive for a random vector $\bm X=(X_1, \ldots, X_n)$ if
\begin{equation}
 \label{eq-250922}
   \VaR_\alpha \(\sum^n_{i=1} X_i\) \le \sum^n_{i=1} \VaR_\alpha(X_i),\quad \alpha \in (0,1).
\end{equation}
If the inequality in \eqref{eq-250922} is reversed, we say VaR is superadditive for a random vector $\bm X$. From Theorem \ref{th-H}, we conclude that VaR is superadditive for a vector of iid random variables $X_1, \ldots, X_n$ with a common distribution belonging to $\mathcal{H}^\ast$. Recently, \cite{IK25} proved that, in an atomless probability space $(\Omega, \mathscr{F}, \P)$, VaR is subadditive for a random vector $\bm X$ with each component integrable (unnecessarily identically distributed) if, and only if $\bm X$ is comonotonic. This result also gives a new equivalent characterization for the comonotonicity of a random vector. For the definition of comonotonicity and its properties, see \cite{DDGKV02}.

It is interesting to investigate sufficient conditions under which VaR is superadditive for a positive random vector. It is natural to wonder whether we have
$$
    \VaR_\alpha \(\sum^n_{i=1} X_i\) \ge \sum^n_{i=1} \VaR_\alpha(X_i),\quad \alpha \in (0,1).
$$
if $X_1, \ldots, X_n$ are independent random variables with $X_i\in\mathcal{H}^\ast$ or $\mathcal{V}$ for each $i$.

In what follows, define
$$
   \mathcal{D}_n^+ =\{\bm\th\in\R^n: \th_1 \ge \th_2\ge \cdots \ge \th_n\ge 0\}.
$$
Let $\bm X=(X_1, \ldots, X_n)$ be a vector of iid random variables with a common distribution $F\in\mathcal{H}$. Another question is whether
\begin{equation}
 \label{eq-250928}
    \(\eta_n X_n, \eta_n X_n+\eta_{n-1} X_{n-1},\ldots, \sum^n_{i=1} \eta_i X_i\)
  \le_{\rm st} \(\th_n X_n, \th_n X_n+\th_{n-1} X_{n-1},\ldots, \sum^n_{i=1} \th_i X_i\)
\end{equation}
holds whenever $\bm\th, \bm\eta\in \mathcal{D}_n^+$ such $\bm \th \prec_{\rm m} \bm\eta$. By Lemma 1 in \cite{Ma98}, there exist a finite number of vectors $\bm\th^{(0)}, \bm \th^{(1)}, \ldots, \bm \th^{(m)}$ in $\mathcal{D}_n^+$ such that $\bm\th =\bm \th^{(0)} \prec_{\rm m} \bm \th^{(1)}\prec_{\rm m}\cdots\prec_{\rm m} \bm\th^{(m)}=\bm\eta$, and for each $k\in [m]$, $\bm\th^{(k-1)}$ and $\bm\th^{(k)}$ differ only in two coordinates. Thus, to prove \eqref{eq-250928}, it suffices to prove that, for $0<\eta<\th< 1-\th<1-\eta$,
$$
   \big (\eta X_1, \eta X_1 + (1-\eta) X_2\big ) \le_{\rm st} \big (\th X_1, \th X_1 + (1-\th) X_2 \big ).
$$
These two questions are still under our investigation.

\appendix

\setcounter{table}{0}
\setcounter{figure}{0}
\setcounter{equation}{0}
\setcounter{section}{1}
%


\section*{Appendices: Proofs of the main results in Section \ref{sect-3}}
\label{Appendix-A}

\begin{lemma}
  \label{le-250930}
For any $(x, y, \beta) \in (0, 1)^3$, we have
\begin{equation*}
    (1-xy)^\beta \leq (1 - x)^\beta + (1 - y)^\beta - (1 - x)^\beta (1 - y)^\beta.
\end{equation*}
\end{lemma}

\begin{proof}
Define $u = 1/(1-x) > 1$ and $v = 1/(1-y) > 1$, and consider the function
\begin{equation*}
h(u, v) = (u+v-1)^\beta - u^\beta - v^\beta + 1.
\end{equation*}
We aim to show that $h(u, v)\le 0$ for $u > 1$ and $v > 1$. Observe that the partial derivative with respect to $u$ is
\begin{equation*}
   \frac{\partial h(u, v)}{\partial u} = \beta \[ (u+v-1)^{\beta-1} - u^{\beta-1}\].
\end{equation*}
Since $\beta - 1 < 0$ and $u + v - 1 > u$, we have $(u+v-1)^{\beta-1} < u^{\beta-1}$, which implies $\partial h(u,v)/\partial u < 0$. Thus, $h(u,v)$ is strictly decreasing in $u\in (1,\oo)$ for fixed $v\in (1,\oo)$, and hence $h(u, v) \le h(1, v) = 0$. This completes the proof of the lemma.
\end{proof}

\noindent \underline{\emph{Proof of Proposition \ref{pr-250925}}}.\ \
(i)\ \ See \cite{ALO24} for the proof of $\mathcal{G} \subset \mathcal{H}^\ast$. $\mathcal{G}$ is a proper subset of $\mathcal{H}^\ast$ since $\essinf (F)=0$ for $F\in \mathcal{G}$ while $\essinf (F)$ may be positive from $F\in\mathcal{H}^\ast$.

To prove $\mathcal{H}\subset \mathcal{V}$, choose $F\in\mathcal{H}$. Then $\eta(x):=\Fbar(1/x)$ is concave in $x\in\R_{++}$. Denote $\eta(0)=\lim_{x\downarrow 0} \Fbar (1/x)=0$. Then $\eta(x)$ is concave on $\R_+$, which implies $\eta(y)/y$ is decreasing on $\R_{++}$, that is, $x \Fbar(x)$ is increasing on $\R_+$. So, $F\in \mathcal{V}$, implying $\mathcal{H}\subset \mathcal{V}$. Example \ref{ex-250920} shows that $\mathcal{H}$ is a proper subset of $\mathcal{V}$.

To prove $\mathcal{V}\subset \mathcal{H}^\ast$, choose $F\in\mathcal{V}$. Denote $\ell(x)=x \Fbar(x)$. Since $F\in\mathcal{V}$, we have $\ell(x)$ is increasing in $x\in\R_+$. Thus,
$$
   \Fbar\(\frac{1}{x_1}\)+ \Fbar\(\frac{1}{x_2}\) = x_1 \ell\(\frac {1}{x_1}\) + x_2 \ell\(\frac {1}{x_2}\)
   \ge x_1 \ell\(\frac {1}{x_1+x_2}\) + x_2 \ell\(\frac {1}{x_1+x_2}\) = \Fbar\(\frac{1}{x_1+x_2}\)\
$$
for all $(x_1,x_2)\in \R_{++}^2$. This means $F\in\mathcal{H}^\ast$, implying $\mathcal{V}\subset \mathcal{H}^\ast$. Example \ref{ex-250921} shows that $\mathcal{V}$ is a proper subset of $\mathcal{H}^\ast$.


(ii)\ \ Denote $\Gbar(x)=1-F^\beta(x)$. Assume $F\in\mathcal{H}^\ast$. Note that $\Gbar(1/x)= \psi\circ \Fbar(1/x)$, where $\psi(x)=1-(1-x)^\beta$ is concave on $[0,1]$ and, hence, subadditive. Then $\Gbar(1/x)$ is subadditive in $x\in\R_{++}$, i.e., $G\in \mathcal{H}^\ast$.

Next, assume $F\in\mathcal{V}$, i.e., $x\Fbar(x)$ is increasing in $x\in\R_+$. Note that
\begin{align*}
  x \Gbar(x) & = x \Fbar (x) \cdot \frac {1-F^\beta(x)}{1-F(x)} = x\Fbar (x) \cdot \varphi (F(x)),
\end{align*}
where $\varphi(t)=[1-t^\beta]/(1-t)$. It is easy to see that
$$
   \varphi'(t) \stackrel {\rm sgn}{=} 1+(\beta-1) t^\beta - \beta t^{\beta-1} \stackrel {\rm def}{=} \zeta(t),
$$
and $\zeta'(t)=\beta (\beta-1) t^{\beta-2} (t-1)\le 0$ for $t\in [0,1]$. Since $\varphi'(1)=0$, it follows that $\varphi'(t)\ge 0$ for $t\in [0,1]$, that is, $\varphi(t)$ is increasing in $t\in [0,1]$. Thus, $x \Gbar(x)$ is increasing in $x\in\R_+$, i.e., $G\in \mathcal{V}$. 

(iii)\ \ We only consider the case for $\mathcal{G}$ since the other cases are trivial. Let $F\in\mathcal{G}$, i.e.,
\begin{equation*}
   F\(\frac{1}{x_1+x_2}\) \geq F\(\frac{1}{x_1}\) F\(\frac{1}{x_2}\),\quad (x_1, x_2) \in \R^2_{++}.
\end{equation*}
Then, by Lemma \ref{le-250930},
\begin{align*}
  \Fbar^\beta\(\frac{1}{x_1+x_2}\) &\le \[1 - F\(\frac{1}{x_1}\) F\(\frac{1}{x_2}\)\]^\beta \\
      &\le \[\Fbar\(\frac{1}{x_1}\)\]^\beta + \[\Fbar\(\frac{1}{x_2}\)\]^\beta - \[\Fbar\(\frac{1}{x_1}\)\]^\beta \[\Fbar\(\frac{1}{x_2}\)\]^\beta.
\end{align*}
Denote $G=1-\Fbar^\beta$. We have
\begin{align*}
  G\(\frac{1}{x_1+x_2}\) & \ge 1-\[\Fbar\(\frac{1}{x_1}\)\]^\beta - \[\Fbar\(\frac{1}{x_2}\)\]^\beta + \[\Fbar\(\frac{1}{x_1}\)\]^\beta \[\Fbar\(\frac{1}{x_2}\)\]^\beta \\[3pt]
     &= G\(\frac{1}{x_1}\) G\(\frac{1}{x_2}\).
\end{align*}
This means $G \in \mathcal{G}$.

(iv)\ \ The proof for the case of $\mathcal{V}$ is trivial. Now, assume $F\in\mathcal{H}^\ast$. Since $F\le_{\rm hr} G$, we have
\begin{align*}
   \Gbar\(\frac {1}{x}\) + \Gbar \(\frac {1}{y}\) & =  \Fbar\(\frac {1}{x}\)\cdot \frac { \Gbar(1/x)}{\Fbar(1/x)}  +  \Fbar\(\frac {1}{y}\)\cdot \frac {\Gbar(1/y)}{ \Fbar (1/y)} \\[3pt]
  & \ge\[\Fbar\(\frac {1}{x}\) + \Fbar\(\frac {1}{y}\)\] \frac { \Gbar(1/(x+y))}{ \Fbar(1/(x+y))} \\[3pt]
  & \ge \Fbar\(\frac {1}{x+y}\) \frac { \Gbar(1/(x+y))}{ \Fbar(1/(x+y))}
      = \Gbar\(\frac {1}{x+y}\),\quad (x,y)\in \R_{++},
\end{align*}
implying $G\in\mathcal{H}^\ast$. \hfill $\Box$
\medskip

\noindent \underline{\emph{Proof of Example \ref{ex-250926}}}.\ \
First, we prove $F^\beta \notin \cal H$ for $\beta= 0.5$. Define $\eta_\beta(x)= 1- F^\beta(1/x)$. Then
\begin{align*}
    \eta_\beta''(x) & = \frac {\beta F^{\beta-2}(1/x)}{x^4} \[ (1-\beta) f^2\(\frac {1}{x}\) -2 x f\(\frac {1}{x}\) -f'\(\frac {1}{x}\)\] \\
    & \stackrel {\rm sgn}{=} (1-\beta) f^2\(\frac {1}{x}\) -2 x f\(\frac {1}{x}\) -f'\(\frac {1}{x}\) \\
    & \stackrel {\rm sgn}{=} \(1-\beta -2\pi\log x\) f\(\frac {1}{x}\) -x \\
    & \stackrel {\rm sgn}{=} \frac {1-\beta -2\pi\log x}{\pi [1+(\log x)^2]} -1 \\
    & \stackrel {\rm sgn}{=} 1-\beta -\pi (1+\log x)^2.
\end{align*}
For $\beta = 0.5$, we find $\eta''_{1/2}(e^{-1}) > 0$, which implies $F^{1/2} \notin \cal H$.

Next, we prove $N(x)\le 0$ for all $x\in\R_{++}$, where
$$
     N(x)=\log \Fbar (x)) + \frac{x f(x)}{\Fbar (x)}.
$$
A straightforward calculation yields
\begin{equation*}
    N(x) = \log\left(\frac{1}{2}-\frac{\arctan(\log x)}{\pi}\right) + \frac{1}{\pi \(1+(\log x)^2\) \[1/2-(1/\pi) \arctan(\log x)\]}.
\end{equation*}
Let $y = 1/2 - (1/\pi) \arctan(\log x)$. Then $y \in (0,1)$, and it remains to prove that
\begin{equation*}
    h(y) := \log y + \frac{\sin^2(\pi y)}{\pi y} \leq 0.
\end{equation*}
We now show that $\psi(y) := \pi y h(y) \leq 0$ for $y \in (0,1)$. Note that $\lim_{y \to 0} \psi(y) = \psi(1) = 0$ and
\begin{equation*}
    \psi'(y) = \pi \left[\log y + 1 + \sin(2\pi y)\right].
\end{equation*}
Setting $\psi'(y) = 0$, we find that the equation has three roots: $y_1 \in (0.15, 0.16)$, $y_2 \in (0.56, 0.58)$, and $y_3 \in (0.84, 0.85)$. Furthermore, $\psi'(y) \leq 0$ on $(0, y_1) \cup (y_2, y_3)$ and $\psi'(y) \geq 0$ on $(y_1, y_2) \cup (y_3, 1)$. Therefore, $\psi(y)$ is decreasing on $(0, y_1) \cup (y_2, y_3)$ and increasing on $(y_1, y_2) \cup (y_3, 1)$. Since $\psi(y_2) < 0$, it follows that $\psi(y) \leq 0$ for all $y \in (0,1)$.

\medskip

\noindent \underline{\emph{Proof of Proposition \ref{pr-250928}}}.\ \
(vi)\ \ Assume $X, Y\in\mathcal{V}$, and denote by $H$ the distribution function of $\max\{X, Y\}$. Then $H(x)=F_X(x) F_Y(x)$ for all $x\in\R_+$. Since $x \Fbar_X(x)$ and $x \Fbar_Y(x)$ are increasing in $x\in\R_+$, it follows that
$$
  x \overline{H}(x) = x \[1-\(1-\Fbar_X(x)\) F_Y(x)\] =x \Fbar_Y(x) + x \Fbar_X (x) F_Y(x)
$$
is also increasing in $x\in\R_+$. Thus, $H\in\mathcal {V}$. Next, assume $X, Y\in\mathcal{H}^\ast$. Then,
\begin{align*}
  \overline{H}\(\frac {1}{x}\) + \overline{H}\(\frac {1}{y}\)
    & = \Fbar_X\(\frac {1}{x}\) +\Fbar_X\(\frac {1}{y}\) + \Fbar_Y\(\frac {1}{x}\) F_X\(\frac {1}{x}\)
        + \Fbar_Y\(\frac {1}{y}\) F_X\(\frac {1}{y}\)  \\
    & \ge \Fbar_X\(\frac {1}{x+y}\) + \Fbar_Y\(\frac {1}{x}\) F_X\(\frac {1}{x+y}\)
        + \Fbar_Y\(\frac {1}{y}\) F_X\(\frac {1}{x+y}\) \\
     & \ge \Fbar_X\(\frac {1}{x+y}\) + \Fbar_Y\(\frac {1}{x+y}\) F_X\(\frac {1}{x+y}\)  \\
     & =  \overline{H}\(\frac {1}{x+y}\),\qquad (x,y)\in\R_{++},
\end{align*}
implying $H\in\mathcal{H}^\ast$.

(vii)\ \  Note that the distribution function of $\max\{X-c, 0\}$ is give by $G(x)=F_X(x+c)$ for $x\in\R_+$.  First, assume $F_X\in\mathcal{V}$. Then $ x \Gbar(x)=x \Fbar_X(x+c)=(x+c)\Fbar_X (x+c) -c \Fbar_X(x+c)$ is increasing in $x\in\R_+$, implying $G\in\mathcal{V}$.

Second, assume $F_X\in\mathcal{H}$. It suffices to show that $\Fbar_X (c+1/x)$ is concave in $x\in\R_{++}$. Note that $\Fbar_X (c+1/x)=\eta\circ \tau(x)$, where $\eta(x)=\Fbar_X(1/x)$ and $\tau(x)=x/(1+cx)$. Since both $\eta(x)$ and $\tau(x)$ are increasing concave, it follows that $\Fbar_X (c+1/x)$ is concave in $x\in\R_{++}$.

Third, assume $F_X\in\mathcal{H}^\ast$. To prove $G\in\mathcal{H}^\ast$, it suffices to show that
\begin{equation}
 \label{eq-250911}
  \Fbar_X\(\frac {1}{x}+c\) +\Fbar_X\(\frac {1}{y}+c\) \ge \Fbar_X\(\frac {1}{x+y}+c\)
\end{equation}
for any $(x,y)\in\R_{++}^2$. Denote $x^\ast =x/(1+cx)$, $y^\ast =y/(1+cy)$ and $a^\ast=(x+y)/(1+(x+y)c)$. It is easy to see that $x^\ast + y^\ast \ge a^\ast$. Since $F_X\in\mathcal{H}^\ast$, we have
\begin{align*}
  \Gbar\(\frac {1}{x}\) +\Gbar\(\frac {1}{y}\) & = \Fbar_X\(\frac {1}{x^\ast}\) +\Fbar_X\(\frac {1}{y^\ast}\)
     \ge \Fbar_X\(\frac {1}{x^\ast+y^\ast}\) \\
     & \ge \Fbar_X\(\frac {1}{a^\ast}\) \ge \Gbar\(\frac {1}{x+y}\),
\end{align*}
implying \eqref{eq-250911}. This proves $G\in\mathcal{H}^\ast$.

Fourth, assume $F_X \in \mathcal{G}$. Since $x^\ast + y^\ast \ge a^\ast$, we have
$$
   G\(\frac{1}{x}\) G\(\frac{1}{y}\) = F_X\(\frac {1}{x^\ast}\) F_X\(\frac {1}{y^\ast}\) \le F_X\(\frac {1}{x^\ast+y^\ast}\) \le F_X\(\frac{1}{a^\ast}\) = G\(\frac{1}{x+y}\),
$$
which implies $G\in\mathcal{G}$.

(viii)\ \ Denote by $H$ the distribution function of $(X-Y)_+$. First, assume $F_X\in\mathcal{V}$, which implies $x \Fbar_X(x+z)$ is increasing in $x\in\R_+$ for each $z\in\R_+$. Then
$$
   x \overline{H} (x) = \int^\oo_0  x \Fbar_X (x+z) \d F_Y(z)
$$
is increasing in $x\in\R_+$, implying $H\in\mathcal{V}$.

Second, the proofs for $F_X\in\mathcal{H}$ and $F_X\in \mathcal{H}^\ast$ directly follow from part (vii).

Third, assume $F_X\in\mathcal{G}$. To prove $H\in\mathcal{G}$, it suffices to show that
\begin{equation}
 \label{eq-250901}
  H\(\frac {1}{x}\) H\(\frac {1}{y}\) \le H\(\frac {1}{x+y}\)
\end{equation}
for any $(x,y)\in\R_{++}^2$. For fixed $(x,y)\in \R_{++}^2$, $F_X(z+1/x)$ and $F_X(z+1/y)$ are both increasing in $z\in\R_+$ and, hence, $F_X(Y+1/x)$ and $F_X(Y+1/y)$ are positively associated \citep{EPW67}. Thus,
$$
   \E \[F_X\(\frac {1}{x}+Y\)\]\cdot  \E \[F_X\(\frac {1}{y}+Y\)\] \le  \E \[F_X\(\frac {1}{x}+Y\) F_X\(\frac {1}{y}+Y\)\].
$$
Consequently, we have
\begin{align}
   H\(\frac {1}{x}\) H\(\frac {1}{y}\) &= \int^\oo_0  x \Fbar_X \(\frac {1}{x}+z\) \d F_Y(z)\cdot \int^\oo_0  x \Fbar_X \(\frac {1}{y}+z\) \d F_Y(z) \nonumber \\
   & =\E \[F_X\(\frac {1}{x}+Y\)\]\cdot  \E \[F_X\(\frac {1}{y}+Y\)\] \nonumber \\
   & \le  \E \[F_X\(\frac {1}{x}+Y\) F_X\(\frac {1}{y}+Y\)\] \nonumber \\
   & \le  \E \[F_X\(\frac {1}{x+y}+Y\)\]  \label{eq-250907} \\
   & =H\(\frac {1}{x+y}\), \nonumber
\end{align}
where \eqref{eq-250907} follows from part (vii) for $\mathcal{G}$. This proves \eqref{eq-250901}.

(ix)\ \ For any distribution $F_X$ from one of $\mathcal{V}$, $\mathcal{H}$ and $\mathcal{H}^\ast$, $F_X$ is heavily-tailed and thus $\Fbar_X(c)>0$. Hence, $[X|X>c]\notin\mathcal{G}$ follows since its essential infimum is not zero. The remaining proof is trivial by observing that the $\P(X>x|X>c)=\min\{\Fbar_X(x)/\Fbar_X(c), 1\}$ for $x\in\R_+$. \hfill $\Box$

\section*{Funding}
	
Z. Zou gratefully acknowledges financial support from National Natural Science Foundation of China (No. 12401625), and the Fundamental Research Funds for the Central Universities (No. WK2040000108). T. Hu gratefully acknowledges financial support from the National Natural Science Foundation of China (No. 72332007, 12371476).

\section*{Disclosure statement}

No potential conflict of interest was reported by the authors.


\end{document}